\documentclass[10pt]{article}
\usepackage{amssymb}
\usepackage{amsmath}
\usepackage[amsmath,thmmarks,thref,standard]{ntheorem}
\usepackage[T1]{fontenc}
\usepackage{fancyhdr}
\newtheorem{lem}{Lemma}[section]
\newtheorem{prop}[lem]{Proposition}
\newtheorem{theo}[lem]{Theorem}
\newtheorem{rem}[lem]{Remark}
\newtheorem{defi}[lem]{Definition}

\newtheorem{ass}[lem]{Assumption}
\newtheorem{cor}[lem]{Corollary}

\renewcommand{\epsilon}{\varepsilon}
\newcommand{\R}{\mathbb R}
\newcommand{\N}{\mathbb N}
\newcommand{\C}{\mathbb C}

\newcommand{\Z}{\mathbb Z}
\newcommand{\E}{\mathbb E}
\newcommand{\B}{\mathcal{B}}

\newcommand{\M}{\mathcal{M}}
\newcommand{\eps}{\varepsilon}

\newcommand{\limes}{\lim_{n\rightarrow\infty}}

\begin{document}

\title{Ornstein-Uhlenbeck Equations with time-dependent coefficients\\ 
and Lévy Noise in finite and infinite dimensions}
\author{Florian Knäble\\ Universität Bielefeld\\E-mail: f.knaeble@googlemail.com}
\date{January 2009}
\maketitle

\begin{abstract}
We solve a time-dependent linear SPDE with additive Lévy noise in the mild and weak sense. Existence of a generalized invariant measure for the associated transition semigroup is established and the generator is characterized on the corresponding $L^2$-space. The square field operator is calculated, allowing to derive a Poincaré and a Harnack inequality.  
\end{abstract}
\section{Introduction}

Recently, there seems to be a growing interest in the study of semilinear non-autonomous stochastic evolution equations in infinite dimensional spaces, see e.g. \cite{veraar} and the references therein. In the so-called semigroup approach which we also pursue in this paper, this is of course closely connected to solving the underlying non-autonomous Cauchy problem, covered by the theory of evolution semigroups or evolution families, see e.g. \cite{evo}.
So, let us consider the following equation:
\begin{equation} \label{intro}
\left\{
\begin{array}{rcl}
dX_t& = &  (A(t)X_t + f(t)) dt + B(t)dZ_t \\
X(s)& = & x
\end{array} \right.
  \end{equation}
  on a Hilbert space $H$, where $A_t:H\to H$ are linear operators and all coefficients are T-periodic.\\
From a mathematical point of view it is natural to study \eqref{intro} with an additional Lipschitz non-linearity, because this case can still be covered if one attacks equation \eqref{intro} using contraction methods (e.g. Picard-Lindelöf iteration). Our main interest, however, in the linear case \eqref{intro}, is to exploit a lot of explicit formulae (e.g. for the transition semigroup, which are in fact time-dependent versions of generalized Mehler formulae) to get as precise as possible information about the solution and the corresponding Kolmogorov equations.
%, cf. Remark \ref{mehler} below
Concerning the time-dependent (possibly unbounded) operators $A_t$, for simplicity, we only consider the case of a common domain. The general case is technically harder, but we expect our results to hold true. This will be the subject of further study.

For $H=\R^d$ and $Z_t$ a d-dimensional Brownian motion, \eqref{intro} was studied intensively in \cite{luna}.
 Inspired by their paper, our work consists primarily in generalizing their results to the case $H$ infinite-dimensional and $Z_t$ a Lévy process. A number of our arguments are adapted from \cite{luna}, although the Lévy setting forces us to work more heavily with Fourier transforms and the infinite dimensional setting requires extra care.
Nevertheless, we succeed in proving the following results: existence and uniqueness of the solution, explicit calculation of its characteristic function, proof of the Chapman-Kolmogorov equations, existence and uniqueness of an evolution system of measures, existence and form of the generator (including a result on its spectrum), precise form of the square field operator, proof of a 
Poincaré and a Harnack inequality.
  
In the autonomous case, the semigroups associated to this kind of equation are known as the generalized Mehler semigroups mentioned above and they are already well understood. Invariant measures are established in \cite{mich} and \cite{GMS}, generators are examined in \cite{lescot} and the square field operator is identified in \cite{sqfo}. 
 
The paper is organized as follows.  
In Section \ref{levy} we shortly review the Lévy-Ito decomposition and the Lévy-Khinchine representation.
In Section \ref{solve} we develop the necessary theory of integration to give sense to our solution, whose existence is established in a rather standard way.
Then, in Section \ref{sgim}, we calculate the Fourier transforms of our solutions and thus determine its transition semigroup explicitly, which is a two-parameter semigroup, since the equation is non-autonomous.
In this case the concept of an invariant measure has to be generalized to allow for a whole collection of measures - a so-called \textit{evolution system of measures}- which are invariant in an appropriate sense. We prove the existence of such a system under the condition that the $A_t$ generate an exponentially stable semigroup and provided that the Lévy symbol is sufficiently smooth. In contrast to the Gaussian setting, the construction of limit measures is more delicate in our case.

Then (as usual) we turn the problem into an autonomous one by enlarging the state space, allowing for a one-parameter semigroup. Via the evolution system of measures and thanks to the periodicity of the coefficients we are able to construct a unique invariant measure for this semigroup as in \cite{luna}. Thus we can introduce the $L^2$-space with respect to the invariant measure where the semigroup is then shown to be strongly continuous.

Section \ref{gen} is dedicated to an analysis of the generator.
We identify an explicit domain of uniqueness for the generator and its action on it, thus identifying it as an pseudo-differential operator in infinitely many variables. 
Subsequently, in Section \ref{asy}, we establish the form of its square field operator which yields a generalization to the crucial ``integration by parts formula'' in \cite{luna}. Then we prove an estimate for the square field operator, that allows us to obtain a Poincaré and a Harnack inequality for our semigroup. %This also implies results for the original two-parameter semigroup.
\\
\noindent
\textbf{Acknowledgement}
I wish to thank Prof. Dr. Röckner for introducing me to this interesting subject and for many helpful suggestions.

\section{Lévy Processes}\label{levy}
In the following let be $H$ a separable Hilbert space with scalar product $\langle \cdot,\cdot \rangle:=\langle \cdot,\cdot \rangle_H$ and norm $\| \cdot \|:=\| \cdot \|_H$.  
An $H$-valued stochastic process L adapted to a 
filtration $(\mathcal{F}_{t})_{t\geq 0}$ is called a \textsl{L\'{e}vy process} if it has independent and stationary increments, is stochastically continuous, and we have $P(L=0)=1$.
We say that $A \in \B(H)$ is bounded below if $ 0 \notin \bar{A}$\newline
We denote by $N(t,A)$ the (random) number of "jumps of size $A$" up to time $t$, that is 
$N(t,A):= \mbox{card}\{ 0\leq s \leq t | \Delta L_s \in A  \}$ \\
If $A$ is bounded below, then $N(t,A)$ is a Poisson process, with intensity $M(A)$, where $M(A):=\E[N(1,A)]$.\\
\begin{defi}
A measure $M$ on $H$ with :
$$\int_{H\backslash\{0\}} \min(1,\|x\|^2) M(dx) <\infty$$
is called a Lévy measure.
\end{defi}

\begin{theo}[Lévy-Ito Decomposition]\label{lito}
If $L$ is an $H$-valued Lévy process, there is a drift vector $b \in H$, a $R$-Wiener process $W_R$ on $H$, such that $W_R$ is independent of $N_t(A)$ for any $A$ that is bounded below and we have:
\[L_t=bt+W_R(t)+\int_{\|x\|<1}x(N_t(dx)-tM(dx))+\int_{\|x\|\geq1}xN_t(dx)\]
where $N_t$ is the Poisson random measure associated to $L$, and $M$ the corresponding Lévy measure. 
\end{theo}
\begin{proof} See e.g. \cite{alb} Theorem 4.1
\end{proof}

\begin{theo}[Lévy-Khinchine Representation]\label{kin}
If $L$ is an $H$-valued Lévy process with Lévy-Ito decomposition as in \ref{lito}, then its characteristic function takes the form: $\E[e^{i\langle L_t,u \rangle}]= e^{t\lambda(u)}$ and
\begin{equation}\label{triple}
\lambda(u)= i\langle b,u \rangle -\frac{1}{2}\langle u,Ru\rangle + \int_{H/\{0\}}\left[e^{i\langle u,x\rangle} -1 -i\langle u,x\rangle \chi_{\{\|x\|\leq 1\}}\right] M(dx) 
\end{equation}
\end{theo}
\begin{proof}
See \cite{linde} Theorem 5.7.3
\end{proof}

Since a measure is characterized by its Fourier transform we will say that a measure $\mu$ is associated to a triple $[b,R,M]$ if its characteristic exponent has the form \eqref{triple}.

\begin{rem}
Actually the Lévy-Khinchine representation holds not only for Lévy processes but for any infinitely divisible random variable. (See \cite{sato} for an account of infinite divisibility) Moreover, Lévy processes and infinitely divisible measures can be brought in a one to one correspondence. 
In particular the converse of \ref{kin} is true: any function of the form
$$\exp\left\{i\langle b,u \rangle -\frac{1}{2}\langle u,Ru\rangle + \int_{H/\{0\}}\left[e^{i\langle u,x\rangle} -1 -i\langle u,x\rangle \chi_{\{\|x\|\leq 1\}}\right] M(dx) \right\}$$
is the characteristic function of a measure.
 \end{rem}

\section{Solving the generalized Ornstein-Uhlenbeck equation}\label{solve}
\subsection[Stochastic Integration]{Stochastic Integration with respect to Lévy martingale measures}
In this subsection we follow \cite{apple}, where the proofs of all results can be found. 
Let be $\mathcal R_1$ the ring of all Borel subsets of the unit ball of $H$ which are bounded below.
\begin{defi}
A Lévy martingale measure on a Hilbert space $H$ is a set function $M:\R_+\times \mathcal R_1 \times \Omega \to H$ satisfying:
\begin{itemize}
\item $M(0,A)=0$ almost surely for all $A\in \mathcal R_1$
\item $M(t,\emptyset)=0$ almost surely
\item almost surely we have:$M(t,A \cup B)=M(t,A)+M(t,B)$ for all $t$ and all disjoint $A$, $B \in \mathcal R_1$
\item $M(t,A)_{\{t\geq 0\}}$ is a square-integrable martingale for each $A\in \mathcal R_1$
\item  if $ A\cap B=\emptyset \;M(t,A)_{\{t\geq 0\}}$ and $M(t,B)_{\{t\geq 0\}}$ are orthogonal, that is:
$\left\langle M(t,A),M(t,B)\right\rangle$ is a real-valued martingale for every $A,B \in \mathcal R_1$
\item $\sup \{\E[\|M(t,A)\|^2] \,|\,A\in\mathcal B(S_n)\}<\infty \qquad \mbox{for every}\quad n\in\N$
\item for every sequence $A_j$ decreasing to the empty set such that $A_j \subset \mathcal B(S_n)$ for all $j$ we have:
    $\qquad\lim_{ j\to \infty} \E[\|M(t,A_j)\|^2] < \infty$
\item for every $s<t$ and every $A\in \mathcal R_1$ we have that $M(t,A)-M(s,A)$ is independent of $\mathcal F_s$
\end{itemize}
\end{defi} 
\begin{prop}
$M(t,A)=\int_A x \tilde{N}_t(dx)$ is a Lévy martingale measure on $H$ for every $A \in \mathcal R_1$.
\end{prop}
%\begin{proof}
%see \cite{stefan} Theorem 2.5.2
%\end{proof}
Similarly as a Wiener process is characterized by its covariance operator, we can describe the covariance structure of a Lévy martingale measure by a family of operators parametrized by our ring $\mathcal R_1$.
\begin{prop} \label{nuclear}
\[\E[|\langle M(t,A),v\rangle|^2] = t \langle v,T_A v\rangle\]
for all $t\geq 0, \;v\in H \;A\in \mathcal R_1 $, where the operators $T_A$ are given by \\
$T_A v:= \int_A T_xv \nu(dx)$ and $T_xv:= \langle x,v\rangle x$.
\end{prop}
%\begin{proof}
%see \cite{stefan} Theorem 2.5.4
%\end{proof}
We will establish only a limited theory of integration, as for our purposes it will be sufficient to integrate deterministic operator valued functions.
We do not even need them to depend on the jump size.
The procedure is the same as for Brownian motion,
so let us introduce the space of our integrands, the approximating simple functions, and state how the integral is defined for them.
For convenience, we set $M([s,t],A):= M(t,A)-M(s,A)$.
\begin{defi}
Let $H'$ be another real separable Hilbert space.\\
Let $\mathcal H^2:=\mathcal H^2(T_-,T_+)$ be the space of all $R:[T_-,T_+]\to \mathcal L(H,U)$ such that $R$ is strongly measurable and we have:
$$ \|R\|_{\mathcal H^2}:= \left(\int_{T_-}^{T_+}\int_{\|x\|<1} tr(R(t)T_xR^*(t))\nu(dx)ds\right)^\frac{1}{2}<\infty$$  
\\[5mm]
\noindent
Let $\mathcal S$ be the space of all $R\in \mathcal H^2$ such that 
$$ R=\sum_{i=0}^n R_i\; \chi_{(t_i,t_{i+1}]} \chi_A $$
where $T_-=t_0<t_1<...<t_{n+1}=T_+$ for some $n\in \N$, where each $R_i\in \mathcal L(H,H')$ 
and where $A \in \mathcal R$\\[5mm]
\noindent
For each $R \in \mathcal S,\quad t\in[T_+,T_-]$ define the stochastic integral as follows:
$$ I_t(R):=\sum_{i=0}^n R_i M([t_i\wedge t,t_{i+1}\wedge t],A)$$
\end{defi}

\begin{prop}
The space $\mathcal H^2$ with inner product 
$$ \langle R,U \rangle := \int_{T_-}^{T_+}\int_{\|x\|<1} tr(R(t)T_xU^*(t))\nu(dx)ds$$
is a Hilbert space.
\end{prop}
%\begin{proof}
%see \cite{harry} Lemma 1.2
%\end{proof}
\begin{prop}\label{dense}
The space $\mathcal S$ is dense in $\mathcal H^2$.
\end{prop}
%\begin{proof} see \cite{harry} Lemma 1.3
%\end{proof} 
\begin{prop} \label{isometry}
We have for any $R \in \mathcal S : \E[I_t(R)]=0$ and 
$$\E[\|I_t(R)\|^2]=\int_{T_-}^{t} \int_A \mbox{tr}(R(s)T_xR^*(s))\nu(dx)ds=\|\chi_{[T_-,t]}R(t)\|_{\mathcal H^2}^2$$
So for $t$ fixed, $\;I_t:\mathcal S \to L^2(\Omega,\mathcal F,P;H)$ is an isometry. 
\end{prop}
%\begin{proof}
%\end{proof}
So we can isometrically extend the operator $I_t$ from $\mathcal S$ to its closure $\mathcal H^2$. 
\subsection{Stochastic Convolution}
We want to give meaning to the integral 
\[X_{U,B}:=\int_s^tU(t,r)B(r)dL(r)\]
which we will call a stochastic convolution. Here $L$ is a $H$-valued Lévy process and we have $U(t,r)\in\mathcal L(H), B(r)\in\mathcal L(H)\; \forall \;s\leq r\leq t$.
In anticipation of the assumptions in section \ref{sgim} we will pose the following conditions:
\begin{itemize}
\item $\sup_{r\in\R} \|B(r)\|_{\mathcal L(H)}<\infty$ 
\item  there is $M>0,\omega>0$ such that : $\|U(t,r)\|_{\mathcal L(H)}<Me^{-\omega(t-r)}$
\item $r\mapsto B(r) $ is measurable and $r \mapsto U(t,r)$ is measurable for any fixed $t$
\end{itemize} 
\begin{prop}
If $U$ and $B$ are as above, the stochastic convolution exists in the following sense:
%We write, according to the Lévy-Ito decomposition \ref{lito}:
\begin{equation}\begin{split}\label{conv}
&\int_s^tU(t,r)B(r)dL(r)\\
&=\int_s^tU(t,r)B(r)b\; dr +\int_s^t\int_{\|x\|\geq1}U(t,r)B(r)x \;N_r(dx)\\
&+\int_s^tU(t,r)B(r)dW_Q(r)+\int_s^t\int_{\|x\|<1}U(t,r)B(r)x\;\tilde{N}_r(dx)\end{split}\end{equation}
\end{prop}
\begin{proof}
The proof is analogous to the one of Theorem 6 in \cite{apple} where $U(t,s)=S(t-s)$ for a strongly continuous semigroup $S$: 
The first term in \eqref{conv} is well defined as a simple Bochner integral, and the second as a finite random sum.
For the other two terms, it is straightforward to check under our assumptions, that the integrands are such that the respective isometries apply. 
\end{proof}

\subsection{Existence of the Mild Solution}
In the following we will have to deal with a non-autonomous abstract Cauchy problem - non-autonomous means we are not in the framework of strongly continuous semigroups anymore. This implies in particular, that we have no easy characterization of well-posedness in the sense of the Hille-Yosida theorem available. There are different, yet technical, approaches (see \cite{neid} and the references therein for a recent overview), but since this subject is not in the primary interest of our work, we assume that the problem is well posed. This is closely related to the notion of \textit{evolution semigroups}. Our definition is taken from \cite{ACP} \\

We consider the following non-autonomous generalisation of the Langevin equation:
\begin{equation} \label{sde}
\left\{
\begin{array}{rcl}
dX_t& = &  (A(t)X_t + f(t)) dt + B(t)dL_t \\
X(s)& = & x
\end{array} \right.
  \end{equation}
where $  B:\R \rightarrow \mathcal{L}(H)$ is strongly continuous and bounded in operator norm, $f: \R \rightarrow H$ is continuous,
$L(t)$ is an $H$-valued Lévy-process and where the $A(t)$ are linear operators on $H$ with common domain $D(A)$ and \\ $A: \R\times D(A) \rightarrow H$ is such that
 we can solve the associated non-autonomous abstract Cauchy problem 
\begin{equation} \label{cauchy}
\left\{
\begin{array}{rcl}
dX_t& = &  (A(t)X_t + f(t)) dt \\
X(s)& = & x
\end{array} \right.
  \end{equation}
according to the following definitions:
\begin{defi}
An exponentially bounded evolution family on $H$ is a two parameter family $\{U(t,s)\}_{t\geq s}$ of bounded linear operators on $H$ such that we have:
\begin{itemize}
\item[(i)]
$U(s,s)=Id\quad$ and 
 $\quad U(t,s)U(s,r)=U(t,r)\quad$ whenever $r\leq s\leq t$
\item[(ii)] for each $x\in H$,  $\quad(t,s)\mapsto U(t,s)x\quad$ is continuous on $\; s \leq t$
\item[(iii)] there is $M>0$ and $\omega > 0$ such that :
$\; \|U(t,s)\| \leq Me^{-\omega(t-s)} \;,\; s\leq t$
\end{itemize}
\end{defi}

\begin{ass}\label{ass}
There is a unique solution to \eqref{cauchy} given by an exponentially bounded evolution family $U(t,s)$ so that the solution takes the form:
\[ X_t=U(t,s)x +\int_s^t U(t,r)f(r)dr\]  
Moreover, we assume that :
$$ \frac{d}{dt}U(t,s)x=A(t)U(t,s)x$$
\end{ass}

\begin{rem}
Note that in the finite dimensional case, where each $A_t$ is automatically bounded we get the existence of an evolution family that solves \eqref{cauchy}, under the reasonable assumption that $t \mapsto A_t$ is continuous and bounded in the operator norm, by solving the following matrix-valued ODE:
\[ \left\{ \begin{array}{l} \frac{\partial}{\partial t} U(t,s) = A(t) U(t,s)\\ U(s,s)= Id  \end{array} \right.  \]
Existence and uniqueness are assured since $(t,M)\mapsto A(t)M$ is globally Lipschitz in $M$.
This result even holds in infinite dimensions, see \cite{krein}.
\end{rem}

\begin{defi}\label{mild}
Given assumption \ref{ass} we call the process:
\[ X(t,s,x)= U(t,s)x + \int_s^t U(t,r)f(r)dr + \int_s^t U(t,r)B(r) dL_r \]
 a mild solution for \eqref{sde}. 
\end{defi}

\subsection{Existence of the Weak Solution}

We have called the above expression a mild solution, though there is no obvious relation to the equation yet.
Now, we will show that our candidate solution actually solves our equation in the weak sense. The following definition makes this precise, but first we need to strengthen our assumption concerning the common domain of the $A(t)$ a little:
\begin{ass}
We require that the adjoint operators $A^*(t)$ also have a common domain independent of $t$ which we will denote by $D(A^*)$. 
Furthermore, we assume that $D(A^*)$ is dense in $H$ and that we have:
\[\frac{d}{dt}U^*(t,s)y = U^*(t,s)A_t^*y \]
for every $y\in D(A^*)$.
\end{ass}

\begin{defi}
An $H$-valued process $X_t$ is called a weak solution for \eqref{sde} if for every $y\in D(A^*)$ we have:
\begin{equation}\label{defweak}
\langle X_t,y\rangle = \langle x,y\rangle + \int_s^t \langle X_r,A_r^*y\rangle dr + \int_s^t \langle f(r),y\rangle dr +  \int_s^t B^*(r)y    dL_r   
\end{equation}
Here $(B^*(r)y)(h):=\langle B^*(r)y,h\rangle $ so that $B^*(r)y \in  \mathcal L(H, \R)$ and the integral is well defined, since\\ 
$\|B^*y\|^2_{\mathcal H^2}\leq (t-s) \sup_r \|B(r)\|^2 \|y\|^2\sum_k\int_{\|x\|<1} \|T_x^{\frac{1}{2}}e_k\|^2 \nu(dx)<\infty$.  
\end{defi}   

\begin{theo}\label{weaksol}
The mild solution $X_t$ from definition \ref{mild} is also a weak solution for \eqref{sde}. Moreover, it is the only weak solution.
\end{theo}
\begin{proof}
%By the expression for $X_t$ (that already contains an integral) we will have to establish the following equality:
%\begin{multline}\label{weak}
%\left\langle U(t,s)x + \int_s^t U(t,r)f(r)dr + \int_s^t U(t,r)B(r) dL_r,y\right\rangle \\
%= \langle x,y\rangle + \int_s^t \left\langle U(r,s)x + \int_s^r U(r,u)f(u)du + \int_s^r U(r,u)B(u) dL_u,A_r^*y\right\rangle dr \\+ \int_s^t \langle f(u),y \rangle du +  \int_s^t B^*(u)y    dL_u   
%\end{multline}
%Therefore, we calculate :
%\[   \int_s^t\langle U(r,s)x,A^*_ry\rangle dr  = \int_s^t\langle x,U^*(r,s)A^*_ry\rangle dr = \langle x,\int_s^t U^*(r,s)A^*_ry dr\rangle \]
 %              \[ = \left\langle x,\int_s^t \frac{d}{dr}U^*(r,s)y dr\right\rangle 
  %              =\left\langle x, [U^*(t,s) U^*(s,s)]y \right\rangle = \langle U(t,s)x -x ,y \rangle \]
%and furthermore:
%\[ \int_s^t \left\langle \int_s^r U(r,u)f(u)du,A^*_r y\right\rangle dr = \int_s^t \int_s^r \left\langle  f(u),U^*(r,u) A^*_r y\right\rangle du dr \]
%\[ = \int_s^t \int_u^t \left\langle  f(u),U^*(r,u) A^*_r y\right\rangle dr du = \int_s^t \int_u^t \left\langle  f(u),\frac{d}{dr}U^*(r,u) y\right\rangle dr du\]  \nopagebreak
%\[ = \int_s^t  \left\langle  f(u),[U^*(t,u)-Id] y  \right\rangle  du = \int_s^t  \left\langle  U(t,u)f(u), y  \right\rangle  du -\int_s^t  \left\langle f(u), y  \right\rangle  du  \]
%\[=\int_s^t  \left\langle  f(u),\int_u^t \frac{d}{dr}U^*(r,u) y dr \right\rangle  du \]
%Hence, \eqref{weak} reduces to:
Analogous to \cite{apple} Theorem 7. 
After some relatively straightforward calculations the problem is reduced to proving the following equality:
\begin{equation}
 \left\langle\int_s^t U(t,r)B(r) dL_r,y\right\rangle =
 \int_s^t\left\langle\int_s^r U(r,u)B(u) dL_u,A_r^*y\right\rangle dr  + \int_s^t B^*(u)y    dL_u  
\end{equation}
 and we will do so with the help of two lemmas.
 
\begin{prop}\label{fub}[stochastic Fubini]
Let be  $(M,\mathcal M,\mu)$ a measure space with $\mu$ finite. By $G^2(M)$ denote the
space of all $\mathcal L(H,H')$- valued mappings $R$ on $[s,t]\times M$ such that
$(r,m)\mapsto  R(r,m)y$ is measurable for each $y\in H$ and \\
$\|R\|^2_{G^2(M)}:= \int_s^t \int_M \|R(r,m)T_x^{\frac{1}{2}}\|^2 \nu(dx) \mu(dm) dr <\infty \quad$
Then we have:
\[\int_M \left(\int_s^t R(u,m) dL_u \right)\mu(dm) =  \int_s^t \left(\int_M R(u,m)\mu(dm)\right) dL_u \] 
\end{prop}
\begin{proof}
see \cite{apple} Theorem 5 \\
\end{proof}

\begin{lem}\label{weakstrong} %[weak-strong-compatibility]
Let be $R \in \mathcal H^2$ and $y\in H$. Then we have:
\[\left\langle\int_s^t R(r)dL_r ,y \right\rangle = \int_s^t R^*(r)y dL_r \]
\end{lem}
\begin{proof}
see \cite{apple} Theorem 4
\end{proof}
Now we are able to finish our proof of \ref{weaksol}:
\begin{align*}
&\int_s^t\left\langle\int_s^r U(r,u)B(u) dL_u,A_r^*y\right\rangle dr 
\stackrel{\ref{weakstrong}}{=}  \int_s^t\left(\int_s^r B^*(u)U^*(r,u)A_r^*y \; dL_u\right) dr \\
&\stackrel{\ref{fub}}{=}  \int_s^t\left(\int_u^t B^*(u)U^*(r,u)A_r^*y \; dr \right)dL_u %& \mbox{by stochastic Fubini}\\
=  \int_s^t \left(B^*(u)\int_u^t \frac{d}{dr}U^*(r,u)y \;dr \right)   dL_u  \\
&=  \int_s^tB^*(u) [U^*(t,u)-Id]y  \;   dL_u  
\stackrel{\ref{weakstrong}}{=}  \left\langle \int_s^t   [U(t,u)-Id] B(u)    dL_u , y \right\rangle \\
&=  \left\langle \int_s^t   U(t,u)B(u)    dL_u , y \right\rangle  -   \left\langle \int_s^t B(u)    dL_u , y \right\rangle 
\end{align*}
and that is precisely what we had to show.
\end{proof}

\section{Semigroup and Invariant Measure}\label{sgim}
\begin{ass}
From now on, we assume that there exists $T\!>\!0$ such that the coefficients $A, f$ and $B$ in \eqref{sde} are $T$-periodic.  
\end{ass}

Recall that the weak solution for \eqref{sde} takes the following form:
\[ X(t,s,x)= U(t,s)x + \int_s^t U(t,r)f(r)dr + \int_s^t U(t,r)B(r) dL_r \]

As opposed to the Gaussian case we are no longer able to give an easy representation of the law of $X(t,s,x)$, but we can calculate its Fourier transform:

\begin{lem}[characteristic function]\label{ft}
\begin{multline*}
	 	 \E\left[\exp\left(i\left\langle h,X(t,s,x)\right\rangle\right)\right]= \\
	 \exp\left\{i\left\langle h, U(t,s)x + \int_s^t U(t,r)f(r)dr\right\rangle \right\} \exp \left\{ \int_s^t \lambda(B^*(r)U^*(t,r)h)dr \right\}
\end{multline*}
where $\lambda $ is the Lévy symbol of $L$.
\end{lem}

\begin{Proof}:
Straightforward, by using the isometries to approximate the stochastic integral by a sum, and then using independence of increments and the Lévy-Khinchin formula. Details are included in the appendix.
\end{Proof}

The following lemma is a straightforward generalization of the standard monotone class theorem.\nopagebreak
\begin{lem}[complex monotone classes] \label{complexmc}
Let $\mathcal H$ be a complex vector space of complex-valued bounded functions, that contains the constants and is closed under componentwise monotone convergence. Let $\M \subset \mathcal H$ be closed under multiplication and complex conjugation.
Then, all bounded $\sigma(\M)$- measurable functions belong to $\mathcal H$.
\end{lem}

The last and the next result in combination will be particularly useful:
\begin{lem}\label{multsys}
The functions $\mathcal M:=\{e^{i\langle h,x\rangle} ,h\in H\}$ form a complex multiplicative systems that generates the Borel $\sigma$-algebra of $H$.
\end{lem} 
\begin{proof}
It is obvious that $\mathcal M$ is closed under multiplication and complex conjugation.\\
 To show that indeed  $\sigma(\mathcal M)=\B(H)$ we make use of the following lemma: (see \cite{schwartz} page 108)
\begin{lem}\label{borel}
A countable family of real-valued functions on a Polish space $X$ separating the points of $X$ already generates the Borel-sigma-algebra of $X$. 
\end{lem} 
Our countable family will be $\{f_{n,k}(x):=\sin(\langle \frac{1}{n} e_k,x\rangle)\}_{k,n\,\in\N} \subset \mathcal M $ where $\{e_k\}$ is an orthonormal basis of $H$.  \newline
Since the sine function is injective in a neighborhood of zero, and the functions $\langle \frac{1}{n} e_k,x\rangle)$ separate the points of $H$, so do the $f_{n,k}$.
As real and imaginary parts of functions in $\mathcal M$, it is clear, that the sigma-algebra generated by them is included in $\sigma(\mathcal M)$.
\end{proof}

%HALBGRUPPE
Now we will show that our solution induces a two-parameter semigroup, defined as follows:
\begin{defi}
 Whenever $f:H\to \C$ is measurable and bounded, define
\[ P(s,t)f(x):= \E[f(X(t,s,x))]  \]
$P(s,t)$ will be called the two-parameter semigroup (associated to the solution $X$).
\end{defi}
\begin{lem} \label{hg}
 For $f$ as above, we have the following flow property, i.e. $P(s,t)$ satisfies the Chapman-Kolmogorov equations: 
\[P(r,s)P(s,t)f(x)=P(r,t)f(x)\]
\end{lem}
\begin{proof}
We will show the equality for the functions $f_h(x)=e^{i\langle h,x \rangle}$ and extend it with the help of \ref{complexmc}.
First note, that by \ref{ft} we have 
\[P(s,t)f_h(x)= \exp\left\{i\left\langle h, U(t,s)x \!+\!\!\int_s^t U(t,r)f(r)dr\right\rangle \!+\!\! \int_s^t \lambda(B^*(r)U^*(t,r)h)dr \right\}\]
so that: 
\begin{align*}P(r,s)P(s,t)f_h(x) &=\E[P(s,t)f_h(X(s,r,x)] \\
 &=\E\left[\exp\left\{i\left\langle U^*(t,s)h, X(s,r,x) \right\rangle \right\}\right]\\
  &\quad \times\exp\left\{i\left\langle h, \int_s^t U(t,r)f(r)dr\right\rangle 
         + \!\int_s^t \lambda(B^*(r)U^*(t,r)h)dr \right\}\end{align*}
but again \ref{ft} gives us the Fourier transform of $X(s,r,x)$ this time evaluated at $U^*(t,s)h$:
\begin{multline*}= \exp\left\{i\left\langle U^*(t,s)h, U(s,r)x \right\rangle \right\}  \\
\times \exp\left\{i\left\langle U^*(t,s)h ,\int_r^s U(s,q)f(q)dq\right\rangle \right\}  \exp\left\{i\left\langle h, \int_s^t U(t,r)f(r)dr\right\rangle \right\}\\ \times
\exp \left\{ \int_r^s \lambda(B^*(q)U^*(s,q)U^*(t,s)h)dq \right\}
  \exp \left\{ \int_s^t \lambda(B^*(r)U^*(t,r)h)dr \right\}\end{multline*}
Interchanging $U(t,s)$ with the integral, as it is a bounded operator, making use of the semigroup property of $U$ and $U^*$ and combining the integrals yields the result for exponential $f$. By monotone convergence, it is easy to see that the space of all bounded measurable $f$ for which the flow property holds is a complex monotone vectorspace.
Hence, the proof is complete.
\end{proof}

\begin{rem}
Note that \ref{hg} is equivalent to the Markov property for our solution, but in our case its direct proof seems to be even more difficult.
\end{rem}

%SYSTEME von MAßEN
\subsection{Evolution Systems of Measures}\label{4.1}

Since our equation is non-autonomous we cannot hope for a single invariant measure.
What one can still expect in our setting is a so called \textit{evolution system of measures}, a whole family $ \{\nu_t\}_{t \in \R}$
of measures such that for all $s<t$ and all bounded measurable $f$:
\begin{equation}\int_{\R^n}P(s,t)f(x)\nu_s(dx)= \int_{\R^n}f(x)\nu_t(dx) \label{eq:esm} \end{equation}

To assure the existence of such a system, besides assumption \ref{ass} we will henceforth require the following condition to hold:
\begin{ass}\label{stability}
$\int_{\|x\|>1} \|x\| M(dx) <\infty $
\end{ass}
The following well known lemma will give a useful growth condition for the Lévy symbol that will allow us to construct limit measures. A proof is contained in the appendix. 
\begin{lem}\label{fre}
Every Lévy symbol $\lambda$ with a Lévy measure $M$ satisfying \ref{stability} is Fréchet differentiable.
In particular such a $\lambda$ is locally of linear growth.
\end{lem}

\begin{theo}\label{central}
Assume \ref{stability}. 
Then, with $\lambda$ the Lévy symbol of $L$,  the functions 
   \[ \hat{\nu_t}(h) := \exp\left\{i\left\langle h, \int_{- \infty}^t U(t,r)f(r)dr\right\rangle \right\} \exp \left\{ \int_{- \infty}^t \lambda\{B^*(r)U^*(t,r)h\}dr \right\}\]
   are the Fourier transforms of an evolution system of measures.\newline
   This system is $T$-periodic, that is we have $\nu_{T+t} = \nu_t$ for any $t$. \newline
   Any other $T$-periodic evolution system of measures, coincides with the above. 
\end{theo}

\begin{proof}:
To establish $T$-periodicity, first note, that we have: \newline
 $ U(t,s)=U(t+T,s+T)$ for any $s<t$, which follows easily from its defining differential equation and the assumption that $A$ is $T$-periodic. Hence, we get 
  \[\int_{- \infty}^{t+T} U(t+T,r)f(r)dr=\int_{- \infty}^t U(t+T,r+T)f(r+T)dr=\int_{- \infty}^t U(t,r)f(r)dr\]
  and for the other integral the argument is the same.

We have to assure that the integrals above exist. Since $U$ is stable and $f$ is bounded on all of $\R$ (as it is continuous and periodic) we have:
\[ \int_{- \infty}^t \|U(t,r)f(r)\|dr \leq  \int_{- \infty}^t M e^{-\omega (t-r)}\|f\|_{\infty} dr = \frac{M}{\omega}\|f\|_{\infty}\]

As $\lambda $ is Fréchet differentiable it has locally linear growth, so that with $\lambda(0)=0$ we have $\|\lambda(u)\|\leq C \|u\|$ on the bounded range of the argument for some $C > 0$. So with $\|B^*\|$ bounded  we can treat the second integral as the first:

\begin{equation}\label{asd}    \int_{- \infty}^t \|\lambda\{B^*(r)U^*(t,r)h\}\|dr 
%\leq  C \int_{- \infty}^t \|B^*(r)U^*(t,r)h\| dr \\ 
\leq C \sup_{r}\|B^*(r)\| \frac{M}{\omega}\|h\| < \infty 
\end{equation} 
where we have used, that $\|U^*\|=\|U\|$.%\\[5mm]

To show that these functions are indeed Fourier transforms of measures we can make use of Lévy's continuity theorem in the finite dimensional case.
We have just proven pointwise convergence of the Fourier transforms of $P\circ [X(t,s,x)]^{-1}$, and that the limit function is continuous in $0$ follows easily by dominated convergence. Pointwise convergence under the integral is clear by continuity of $\lambda, U $ and $B$ and a majorizing function is found by looking at \eqref{asd} again. \\
In the infinite dimensional case, however, we cannot apply Lévy's continuity theorem, because we are unable to prove continuity in the Sazonov topology.
% for reasons explained in appendix B. 
For a better readability we postpone the somewhat technical alternative to the end of this proof, formulated as a claim.

In order to see that the respective measures constitute an evolution system of measures we will check \eqref{eq:esm} for exponential functions and then extend the result via monotone classes.
%\newline
So if we take $k(x)= e^{i\langle h,x\rangle}$ in \eqref{eq:esm} we get:
\[  \int_{\R^n}k(x)\nu_t(dx) = \hat{\nu_t}(h) \]
by the very definition of Fourier transformation.
\newline
On the other hand we have by \ref{ft}: 
\[ P(s,t)k(x)= \exp\left\{i\left\langle h,U(t,s)x +\int_{s}^t U(t,r)f(r)dr\right\rangle \!+\! \int_{s}^t \lambda\{B^*(r)U^*(t,r)h\}dr \right\} \]
Using the adjoint of $U$ and the fact that Fourier transformation is only with respect to $x$ we obtain by definition of $\hat{\nu_s}$:
\begin{align*} &\int_{\R^n}P(s,t)k(x)\nu_s(dx)\\
&=\hat{\nu_s}(U^*(t,s)h)\; \exp\left\{i\left\langle h, \int_{s}^t U(t,r)f(r)dr\right\rangle + \int_{s}^t \lambda\{B^*(r)U^*(t,r)h\}dr \right\}\\
&= \exp\left\{i\left\langle U^*(t,s)h, \int_{- \infty}^s U(s,r)f(r)dr\right\rangle + \int_{- \infty}^s \lambda\{B^*(r)U^*(s,r)U^*(t,s)h\}dr \right\} \\ 
 & \qquad \qquad \qquad \times\exp\left\{i\left\langle h, \int_{s}^t U(t,r)f(r)dr\right\rangle + \int_{s}^t \lambda\{B^*(r)U^*(t,r)h\}dr \right\}\\
&= \exp\left\{i\left\langle h, \int_{- \infty}^s U(t,s) U(s,r)f(r)dr\right\rangle + \int_{- \infty}^s \lambda\{B^*(r)U^*(t,r)h\}dr \right\} \\ 
 &\qquad \qquad \qquad\times\exp\left\{i\left\langle h, \int_{s}^t U(t,r)f(r)dr\right\rangle + \int_{s}^t \lambda\{B^*(r)U^*(t,r)h\}dr \right\}\\
%where we used in the left part, that linear continuous operators commute with the integral, and in the right part the corresponding ("twisted") semigroup property for the adjoints. 
%\begin{multline*} = \exp\left\{i\left\langle h, \int_{- \infty}^s U(t,r)f(r)dr\right\rangle    
% + i\left\langle h, \int_{s}^t U(t,r)f(r)dr\right\rangle \right\}  \\
% \exp \left\{ \int_{- \infty}^s \lambda\{B^*(r)U^*(t,r)h\}dr  +  \int_{s}^t \lambda\{B^*(r)U^*(t,r)h\}dr \right\}
%\end{multline*}
 &=\exp\left\{i\left\langle h, \int_{- \infty}^t U(t,r)f(r)dr \right\rangle\right\}
\exp \left\{ \int_{- \infty}^t \lambda\{B^*(r)U^*(t,r)h\}dr \right\}  
\end{align*}
but the last line equals $\hat{\nu_t}(h)$ and that is precisely what we had to show.

To prove the full assertion we have to show that \eqref{eq:esm} not only holds for functions of the form $k_h(x):=e^{i\langle h,x\rangle}$,
but for any bounded measurable function.\\
By \ref{multsys} we can apply \ref{complexmc}, because the bounded and measurable functions for which \eqref{eq:esm} holds, form a complex monotone vector space:\\
for constant functions the equality is trivial and that \eqref{eq:esm} holds for monotone limits is essentially an iterated application of Lévy's theorem about monotone convergence.  %Hence, all bounded measurable functions satisfy \eqref{eq:esm}
Hence, the existence of an evolution system of measures is proved. \\[2mm]
To prove uniqueness, let $\{ \nu_s  \}$ be another $T$-periodic family satisfying \eqref{eq:esm}, then it follows by periodicity:
\begin{multline*} \hat{\nu_s}(h) = \quad \hat{\nu_s}(U^*(s+T,s)h) \\
\times \exp\left\{i\left\langle h, \int_{s}^{s+T} U(s+T,r)f(r)dr\right\rangle + \int_{s}^{s+T} \lambda\{B^*(r)U^*(s+T,r)h\}dr \right\}
\end{multline*}
Using the easy to check relations:
\[\int_s^{s+T}U(s+T,r)f(r)dr=\int_{-\infty}^{s}U(s,r)f(r)dr -U(s+T,s) \int_{-\infty}^{s}U(s,r)f(r)dr \]
\begin{multline*} \int_s^{s+T} \lambda\{ B^*(r)U^*(s+T,r)h\}dr = \\
 \int_{-\infty}^{s} \lambda\{ B^*(r)U^*(s,r)h\}dr -\int_{-\infty}^{s} \lambda\{ B^*(r)U^*(s+T,r)h\}dr   
 \end{multline*}
 
 we get: 
 \begin{align*}
 \hat{\nu_s}(h)= & \;\hat{\nu_s}(U^*(s+T,s)h) \exp\left\{i\left\langle h, \int_{- \infty}^s U(s,r)f(r)dr\right\rangle + \int_{- \infty}^s \lambda\{B^*(r)U^*(s,r)h\}dr \right\} \\
   & \times \exp\left\{i\left\langle h, -U(s+T,s)\int_{- \infty}^s U(s,r)f(r)dr\right\rangle 
   - \int_{- \infty}^s \lambda\{B^*(r)U^*(s+T,r)h\}dr \right\}
\end{align*}
or equivalently:
\begin{gather*}
 \hat{\nu_s}(h) \left[\exp\left\{i\left\langle h, \int_{- \infty}^s U(s,r)f(r)dr\right\rangle \right\} \exp \left\{ \int_{- \infty}^s \lambda\{B^*(r)U^*(s,r)h\}dr \right\}\right]^{-1}\\ 
=\hat{\nu_s}(U^*(s+T,s)h) \left[\exp\left\{i\left\langle U^*(s+T,s)h, \int_{- \infty}^s U(s,r)f(r)dr\right\rangle \right\} \right.\\
 \left.\qquad\qquad\qquad\exp \left\{ \int_{- \infty}^s \lambda\{B^*(r)U^*(s+T,r)h\}dr \right\}\right]^{-1}
\end{gather*}
Finding, that the second line is the first, with $h$ replaced by $U^*(s+T,r)h$ we can iterate, since the relation was valid for all $h \in H$.
As $ \| U^*(s+T,r)\| < 1 $ must hold by our stability assumption, all the factors in the second equation will tend to $1$, so $\hat{\nu_s}$ must indeed have the desired form.\\[2mm]
\textbf{Claim}: $\hat{\nu_t}$ is a characteristic function - Hilbert space case\\
 The general idea is the following. In the Gaussian case, it is known that the limit distributions are Gaussian again. In the same manner we take advantage of the fact that, in the Lévy case, our limit distribution is infinitely divisible.
 We proceed here similarly as in \cite{GMS} chapter 3.
  First of all we show that our distributions $P\circ X(t,s,x)$ are infinitely divisible for any $t>s,x$. By the Lévy-Khinchine representation it is sufficient to prove that their characteristic functions have the form \eqref{triple} for some triple $[b,Q,M]$. Therefore, we calculate:
\begin{align*} &\exp\left\{\int_s^t \lambda(B^*_rU^*_{t,r}h)dr \right\}
\\&
=\exp\left\{\int_s^t i\langle b,B^*_rU^*_{t,r}h\rangle dr-\frac{1}{2}\int_s^t \langle B^*_rU^*_{t,r}h,RB^*_rU^*_{t,r}h\rangle dr\right.\\
&+\left. \int_s^t \left(\int_H e^{i\langle x,B^*_rU^*_{t,r}h\rangle}-1-i\langle x,B^*_rU^*_{t,r}h\rangle\chi_{\{\|x\|\leq1\}} M(dx)\right)dr \right\}\end{align*}
For the jump part we have:
\begin{align}  &\int_H e^{i\langle x,B^*_rU^*_{t,r}h\rangle}-1-i\langle x,B^*_rU^*_{t,r}h\rangle\chi_{\{\|x\|\leq1\}} M(dx) \notag\\
&= \int_H e^{i\langle U_{t,r}B_rx,h\rangle}-1-i\langle U_{t,r}B_r x,h\rangle \chi_{\{\|x\|\leq1\}}M(dx)\notag \\
& + \int_H [ -\chi_{\{\|U_{t,r}B_r x\|\leq1\}} +\chi_{\{\|U_{t,r}B_r x\|\leq1\}}] M(dx)       \notag  \\
&= \int_H e^{i\langle x,h\rangle}-1-i\langle  x,h\rangle \chi_{\{\| x\|\leq1\}} M\circ(U_{t,r}B_r)^{-1}(dx)\label{ff} \\
&-\int_H i\langle  U_{t,r}B_rx,h\rangle[\chi_{\{\|x\|\leq1\}} -\chi_{\{\|U_{t,r}B_r x\|\leq1\}}] M(dx)\label{f}
\end{align}
Note that \eqref{ff} is finite because of: (setting $C:=\|U_{t,r}B_r\|_{\mathcal L(H)}$)
\begin{gather*}\int_H (1\wedge \|x\|^2)M\circ(U_{t,r}B_r)^{-1}(dx)=\int_H (1\wedge \|U_{t,r}B_rx\|^2)M(dx)\\
\leq \int_H (1\wedge C^2\|x\|^2)M(dx)\leq C^2 \int_H (\frac{1}{C^2}\wedge \|x\|^2)M(dx)<\infty \end{gather*}
and only in that way we can argue that \eqref{f} must be finite as well.
Thus, we obtain:
  \begin{align*} &\exp\left\{\int_s^t \lambda(B^*_rU^*_{t,r}h)dr \right\}\\
  &=\exp\left\{ i\left\langle \int_s^t U_{t,r}B_rb \,dr,h  \right\rangle -\frac{1}{2} \left\langle h,\int_s^tU_{t,r}B_rRB^*_rU^*_{t,r}h dr\right\rangle \right.\\
&\quad+ \int_s^t\left( \int_H e^{i\langle x,h\rangle}-1-i\langle  x,h\rangle \chi_{\{\| x\|\leq1\}} M\circ(U_{t,r}B_r)^{-1}(dx) \right)dr \\
&\quad-\left.i\left\langle\int_s^t \int_H   U_{t,r}B_rx[\chi_{\{\|x\|\leq1\}} -\chi_{\{\|U_{t,r}B_r x\|\leq1\}}] M(dx)dr,h \right\rangle
\right\}\end{align*}
so that with:
\begin{itemize} 
\item $b(t,s):=\int_s^t U_{t,r}B_rb dr -\int_s^t \int_H   U_{t,r}B_rx[\chi_{\{\|x\|\leq1\}} -\chi_{\{\|U_{t,r}B_r x\|\leq1\}}] M(dx)dr$
\item $Q(t,s):=\int_s^tU_{t,r}B_rRB^*_rU^*_{t,r}dr$
\item $M_{t,s}(A):=\int_s^tM\circ(U_{t,r}B_r)^{-1}(A)dr$ for $0\notin A$
\end{itemize}
 $\exp\left\{\int_s^t \lambda(B^*_rU^*_{t,r}h)dr \right\}$ is associated to the triple $[b(t,s),Q(t,s),M_{t,s}]$, where $Q(t,s)$ is still symmetric and nonnegative and we have :
\begin{gather*}\mbox{tr} Q(t,s)=\sum_k \langle e_k,Q(t,s)e_k\rangle = \sum_k \int_s^t \|\sqrt{R} B^*_rU^*_{t,r}e_k\|^2 \\
= \int_s^t  \|\sqrt{R} B^*_rU^*_{t,r}\|^2_2\leq \int_s^t  \|\sqrt{R} \|^2_2 \|B^*_rU^*_{t,r}\|^2 < \infty \end{gather*}

and $M_{t,s}$ is a Lévy measure, as we have [since $(1\wedge\|x\|^2)\leq (\|x\|\wedge\|x\|^2)$]:
\begin{align*}&\int_s^t \int_H (1\wedge \|x\|^2)M\circ(U_{t,r}B_r)^{-1}(dx) \leq \int_s^t \int_H (\|x\|\wedge \|x\|^2)M\circ(U_{t,r}B_r)^{-1}(dx) \\
&=\int_s^t\int_H (\|U_{t,r}B_rx\|\wedge \|U_{t,r}B_rx\|^2)M(dx) 
\\&
\leq \int_s^t \|U_{t,r}B_r\|\int_H (\|x\|\wedge \|U_{t,r}B_r\|\, \|x\|^2)M(dx)\\
&\leq \int_s^t \max_{r\leq t} \|B_r\|Me^{-\omega(t-r)} \underbrace{\int_H (\|x\|\wedge \max_{r\leq t} \|U_{t,r}B_r\|\, \|x\|^2)M(dx)}_{< \infty \;\mbox{by assumption}\;\ref{stability}} < \infty
\end{align*}
Moreover, we see that we can let $s \to -\infty$ and $Q(t,-\infty)$ will still be trace class as well as $M_{t,-\infty}$ will still be a Lévy measure, because of the exponential stability of $U$. Since we already know that the Fourier transform as a whole converges, convergence of the first part of $b(t,-\infty)$ (which is obvious) implies convergence of the second part. Hence the limit function is associated to a Lévy triple and thus the characteristic function of an infinitely divisible measure. 
\end{proof}

\begin{rem}
The condition that the Lévy symbol is of linear growth is actually stronger than necessary. To assure the existence of the integral in \eqref{asd} it would be even sufficient to have a very weak estimate of the form $|\lambda(u)|= \mathcal O(\sqrt{\|u\|})$. But we were unable to find any other easy to check conditions to control the growth of a Lévy symbol around the origin. Moreover, in the infinite dimensional case, we have made full use of our assumption.
\end{rem}

\subsection{The Invariant Measure $\nu$ and the Space $L^2_*(\nu)$}

Since we are interested in generators, we have to reduce our equation to the autonomous case, so that we obtain a one-parameter semigroup, that we can relate to a generator.
It will turn out, that we can establish an invariant measure for our new semigroup on the extended state space, using our evolution system of measures.

Reduction of non-autonomous problems is a well-known method in the theory of ordinary differential equations.(see e.g. \cite{krein})
We recall that the basic idea is to enlarge the state space, thus allowing to keep track of the elapsed time.    
The reduced problem then looks:
\begin{equation*}
\left\{
\begin{aligned}
dX(t)&= \{A(y(t))X(t) + f(y(t))\}dt + B(y(t))dL(t)  & \qquad  X(0)&=x \\
dy(t)&=dt                                     & \qquad  y(0)&=s 
\end{aligned}
\right.
\end{equation*}

The one-parameter semigroup is then defined as follows:
   \[ P_{\tau} u(t,x):= P_{t,t+\tau}u(t+ \tau,\cdot)(x) :=(P_{t,t+\tau}u_{t+ \tau})(x) \]
   meaning that we apply the two-parameter semigroup to $u$ as a function of $x$ only.
   That the family $\{P_{\tau}\}_{\tau \in \R}$ is indeed a semigroup, follows, of course, from the semigroup property of 
   $\{P_{s,t}\}_{s<t}$ and is a simple calculation:
   \[(P_{\sigma} (P_{\tau} u)) (t,x)= P_{t,t+\sigma} P_{t+\sigma,t+\sigma+\tau} u (t+\sigma+\tau,\cdot)(x)\]
   \[= P_{t,t+\sigma+\tau} u (t+\sigma+\tau,\cdot)(x)=P_{\tau+\sigma} u(t,x) \]

Starting from our evolution system of measures, we will establish an invariant measure for the one-parameter semigroup. On the respective $L^2$-space the semigroup will then be a contraction.
   
   From now on we will require the following assumption to hold:
   \begin{ass}\label{A*}
   $D(A^*)$ is dense in $H$ and we have $U^*(t,s)D(A^*)\subset D(A^*)$ for all $s\leq t$. 
   Moreover, we have $\frac{d}{dt}U^*(t,s)x=U^*(t,s)A^*(t)x$ for all $x\in D(A^*)$
   \end{ass}
   To obtain our invariant measure we need the following lemma:   
   \begin{lem}
   The function $F: (t,A)\mapsto \nu_t(A) \quad t\in \R, \; A\in \mathcal B(H)\;$ is a kernel.
   \end{lem}
   \begin{proof}
   It is not hard to see, that $t\mapsto \nu_t$ is weakly continuous. A monotone class argument then gives the result.
   \end{proof}
    Now, we introduce the space, on which the semigroup will be strongly continuous and the subspace that will be the core for the generator of our semigroup. 
   \begin{defi}
   As $F$ is a kernel from $H$ to $\R$ we can form $\nu:= F \otimes \frac{1}{T}dt$,  a measure on $H \times \R$.
   \begin{align*}L^2_*(\nu):= & \{ f:\R \times H \to \R \;\mbox{measurable}\, |\, f(t+T,x)=f(t,x)\quad \nu -a.e.\\ & \int_{[0,T]\times H} \|f\|^2(y) \nu(dy) <\infty \} \end{align*}
   % Definition der Räume
   \vspace{-5mm}
   \begin{align*} \M:= \mbox{span}_{\C}& \{ f:\R \times H \to \C  \quad |\quad f= \Phi(t)e^{i\langle x,h(t)\rangle}, \;      \mbox{where}\\ & \Phi \in C^1(\R,\R)\,\mbox{and T-periodic} ,\\
    &h \in C^1(\R,H)\, \mbox{ and T-periodic such that}\; \mbox{Im}\;h \subset D(A^*)\} 
 \end{align*}
   $\qquad \; K:=\{\Re(f) \;| \quad f\in \M\}$\\[2mm]
   That is, $K$ comprises the real parts of the functions in $\M$.
    
   \end{defi}
  
  \begin{rem}
   It is not hard to see that $L^2_*$ is a Hilbert space. %Because of the periodicity it is clear, that $(\int_{[0,T]\times H} \|f\|^2(y) \nu(dy))^{\frac{1}{2}}$ is a norm(where we introduce $\nu \;a.e.$-equivalence classes as usual). Given a Cauchy-sequence $f_n$ we consider $f_n^z$ the restriction of $f_n$ to the interval $I_z:=[zT,(z+1)T], z\in\Z$. By Riesz-Fischer we obtain a limit $f^0$ of $f_n^0$ on $[0,T]$, and because of periodicity it is clear that the other restrictions form the same Cauchy sequences, that is: $\lim f_n^z=f^z=f^0 \; \forall z \in \Z$. Hence, the limit function is periodic, and the space is complete.       
   \end{rem}
  
    \begin{lem}\label{Kdense}
    $K$ is dense in $L^2_*(\nu)$.
    \end{lem}
    
    \begin{proof}
    By a monotone class argument. For the convenience of the reader, details can be found in the appendix. 
    \end{proof}
\begin{rem}
Note that to prove the density of $K$, $h\equiv const$ is sufficient, but we will need the $t$-dependence of $h$ later to show that $K$ is $P_{\tau}$-invariant.
\end{rem}   
   \begin{prop} \label{invext} 
   The measure $\nu$ is the unique invariant measure for the semigroup $P_{\tau}$ on $L^2_*(\nu)$.
   The semigroup $P_{\tau}$ is a contraction on $L^2_*(\nu)$.
   \end{prop}
   \begin{proof}
 %  We will prove both invariance and contractivity for functions from $K$ first. Since $K$ is dense, it is then easy to see that they hold on the whole of $L^2_*(\nu)$. Note, however, that we first need invariance to prove contractivity, then extend contractivity by density, and then obtain global invariance by density and contractivity.
   The proof is analogous to the one in \cite{luna}.
   For the convenience of the reader we include it in the appendix.
  \end{proof} 
 %  \begin{rem}
  % Note that we cannot abandon the T-periodicity, because we need translation invariance in our proof above. But as we want a %probability measure, we cannot take Lebesgue measure on the whole of $\R$. The only alternative known to us, is to restrict ourselves to a finite interval by periodicity.
  % \end{rem}
   %\newpage
 \section{Generator and Domain of Uniqueness} \label{gen}
   In this section we prove that the generator is given by a pseudo-differential operator. Compared to the Gaussian case, we have an additional nonlocal part. 
   However, we still obtain a result on the spectrum of the generator, exactly as in the Gaussian case. 
   
   \begin{prop}[strong continuity]
   $P_{\tau}$ is strongly continuous on $L^2_*(\nu)$  
   \end{prop}
   \begin{proof}
   The invariance of $\nu$ and the density of $K$ allow us to use proposition 4.3 from \cite{diri}. Hence, it will be sufficient to show, that $P_{\tau}u \stackrel{t\rightarrow 0}{\longrightarrow} u $ $\nu - a.e.$ For $u(t,x)=\Phi(t)e^{i\langle x,h(t)\rangle}$ we have by \ref{ft}: 
   \begin{multline} \label{K}
   (P_{\tau}u)(t,x) =  \exp \left\{ \int_t^{t+\tau} \lambda(B^*(r)U^*(t+\tau,r)h(t+\tau))dr \right\}\\
  \times \quad \Phi(t+\tau)\exp\left\{i\left\langle h(t+\tau), U(t+\tau,t)x + \int_t^{t+\tau} U(t+\tau,r)f(r)dr\right\rangle \right\} 
\end{multline}
Recalling that $\Phi,h$ and $U$ are continuous and that $U(t,t)=Id$ we obtain the result, since all the integrals vanish.
Note, that by linearity, this extends to general $u\in K$. 
   \end{proof}
   
   \begin{lem}\label{core}
Let $S(t)$ be a strongly continuous semigroup with generator $G$ on a Banach space $B$. Let $C$ be a dense subspace of $D(G)$ that is invariant under $S$. Then $C$ is a core.
\end{lem}
\begin{proof}
see \cite{arendt} page 47
\end{proof}   
   \begin{lem}
   Let assumption \ref{stability} hold.
  Then $K$ is a core for $G$.
   \end{lem}
   \begin{proof}
    Looking closely at \eqref{K} again, one notes that $(P_{\tau}u)(t,x)$ is again of the form $\Psi(t)e^{i\langle x,k(t)\rangle}$ with $\Psi$ and $k$ as follows: \\[1mm]
    $\Psi(t):= \Phi(t+\tau)\exp\left\{i\left\langle h(t+\tau), \int_t^{t+\tau} U(t+\tau,r)f(r)dr\right\rangle \right\}$ \\ $\mbox{}\quad\qquad \times \quad \exp \left\{ \int_t^{t+\tau} \lambda(B^*(r)U^*(t+\tau,r)h(t+\tau))dr \right\}$ \\[1mm]
  $k(t):= U^*(t+\tau,t)h(t+\tau)$ \\[1mm]
    Since by assumption \ref{A*} $k:\R\to D(A^*)$,  $K$ is invariant under $P_{\tau}$, since by linearity it is sufficient to check the invariance for special $u$. Furthermore, we have again by \eqref{K}:
    \begin{multline}\label{generator}
    Gu=\left.\frac{d}{d\tau}P_{\tau}u\right|_{\tau=0}=[\Phi'(t) + i\Phi(t)\langle x,h'(t)\rangle] \,e^{i\langle x,h(t)\rangle} \\
                                            + i\langle x + f(t),A^*(t)h(t)\rangle\, \Phi(t)\,e^{i\langle x,h(t)\rangle}\\
                                            + \lambda[B^*(t)h(t)]\, \Phi(t)\,e^{i\langle x,h(t)\rangle}
    \end{multline}
    Note that we have used the differentiability of $\lambda$.
   
    Seeing, that $K \subset D(G)$ we can apply \ref{core} to prove the assertion. 
   \end{proof}
   Remember, that $(b,R,M)$ is the triple of our Lévy process (see \ref{triple}).
   Set $\Sigma=\sqrt{R}$.\\
   %To obtain a realization of $G$ let us define the following operator on $K$:
\begin{defi}
For $u\in K$ we set:
\begin{multline*}
Lu(t,x):= u_t(t,x) + \langle A(t)x+f(t), \nabla_x u(t,x) \rangle \\
        +\langle B(t)b, \nabla_x u(t,x) \rangle +\frac{1}{2} \mbox{Tr}\{\Sigma^* B^* \nabla_{xx}u(t,x) B \Sigma  \}  \\
        + \int_{\R^d} \{ u(t,x+B(t)y) - u(t,x) - \langle B(t)y,\nabla_x u(t,x)\chi_{\|y\|\leq 1} \}\nu(dy)	
\end{multline*}
where $\nabla_x\;u$ denotes the gradient of $u$, $\; \nabla_{xx}\;u$ denotes the generalized Hessian of $u$ and Tr denotes the trace of an operator.\\
Note that $\nabla_{xx}\;u$ is Hilbert-Schmidt, (see e.g. \cite{kolm} 1.2.4.) so the trace is well defined.
\end{defi}

\begin{lem}[realization of G]\label{realization}
 $\quad L=G_{|K}\;$, so that  $G=\bar{L}$
 \end{lem}
\begin{proof}
Again, we will only check this for special $u$, since we deal with linear operators.
Note, for the calculations involved, that for $u(t,x)=\Phi(t) e^{i\langle x,h(t) \rangle }$ we have $\nabla_x u(t,x) = ih(t) u(t,x) \quad \nabla_{xx} u(t,x) = -h(t) h^*(t) u(t,x) \quad$ and that \\$ \mbox{Tr}(Au u^*A^*) = \langle Au,Au \rangle $. Hence the Lévy-Khinchine formula yields for $\lambda$:
\begin{align*} &u(t,x) \; \lambda(B^*(t)h(t))  \\
& = i\langle B(t)b,h(t) \rangle \; u(t,x) - \frac{1}{2}\langle \Sigma^* B^*(t)h(t), \Sigma^* B^*(t)h(t)\rangle \; u(t,x) \\
&  \quad+ u(t,x) \int_{\R^d} \{e^{i\langle B^*(t)h(t),y\rangle} -1 -i\langle B^*(t)h(t),y\rangle \chi_{\|y\|\leq 1} \} \nu(dy)\\
& =  \langle B(t)b, \nabla_x u(t,x) \rangle -\frac{1}{2} \mbox{Tr} (\Sigma^* B^*(t)h(t)h^*(t)u(t,x)B(t)\Sigma)\\
& \quad+ \int_{\R^d} \{ \Phi(t)e^{i\langle h(t),x\rangle}e^{i\langle h(t),B(t)y\rangle} - u - i\langle h(t)u,B(t)y\rangle \chi_{\|y\|\leq 1} \}\nu(dy)\\
& =  \langle B(t)b, \nabla_x u(t,x) \rangle +\frac{1}{2} \mbox{Tr}\{\Sigma^* B^* \nabla_{xx}u(t,x) B \Sigma  \}\\
&  \quad+ \int_{\R^d} \{ u(t,x+B(t)y) - u(t,x) - \langle B(t)y,\nabla_x u(t,x)\chi_{\|y\|\leq 1} \}\nu(dy)
\end{align*}
That the first two summands in both expressions also coincide is very easy to see.
\end{proof}   

\begin{lem}\label{diss}
For all $u \in D(L)$ we have 
\begin{equation} \int_{[0,T]\times \R^d} Lu(t,x) \nu(dt,dx)=0 \end{equation}
 \end{lem}

\begin{proof}
We will prove this for $u\in K$ first.
As $\nu$ is $P_{\tau}$ invariant, let us consider the equality:
$$ \int_{[0,T]\times \R^d} P_{\tau}u(t,x) \nu(dt,dx)=\int_{[0,T]\times \R^d} u(t,x) \nu(dt,dx)$$
Differentiating both sides with respect to $\tau$ we obtain the result, if we can show that we can interchange integral and differential on the left hand side. We know, that for a fixed $u \in K \quad \tau \mapsto P_{\tau}u(t,x) $ is differentiable $\nu-$ almost everywhere, since $\limes nP_{\frac{1}{n}}u(t,x)=Lu(t,x) \;$in $ L^2_*(\nu)$ implies pointwise $a.e.$-convergence along a subsequence, that we can choose without loss of generality.
Furthermore, for a fixed $u$ the derivative in $\tau$ is given by  $ P_{\tau}Lu(t,x)$, but since every $P_{\tau}$ is a contraction with respect to the supremum norm, and every $Gu \in K$ is bounded we have a uniform bound for the derivatives. Hence we can apply Lebesgue's dominated convergence theorem.\\
Since $K$ is a core, for $u \in D(G)$  we can find a sequence $u_n$ such that $u_n\to u$ in $L^2_*$ and $Gu_n\to Gu$ in $L^2_*$. Taking the limit, we obtain the equality for $u$.   
\end{proof}

Exactly as in \cite{luna} we obtain the following result on the spectrum of $G$:
\begin{cor}
For any $z \in \sigma(G)$ and $k\in \Z$ we have $z+2\frac{\pi}{T}ki \in \sigma(G)$. Moreover 0 is a simple eigenvalue of $G$. 
\end{cor}
\begin{proof}Analogous to \cite{luna} Corollary 5.5
\end{proof}
%For fixed $k \in\Z$ consider the operator $T_k u(t,x)=e^{2k\frac{\pi}{T}it}u(t,x)$. Since $T$ is unitary the spectrum of $G$ is equal to the spectrum of $T_k^{-1}GT_k=G+(2ki\frac{\pi}{T})Id$ where the equality holds, because the factors cancel out everywhere, except for the derivative with respect to $t$ where the product rule applies. This proves the first statement. \\
%Since every unique invariant measure is ergodic, we also have the equivalent property (see \cite{erg}):\\
%If $u\in L^2_*$ fulfills $P_{\tau}u=u$ for every $\tau>0$ then $u$ is equal to a constant in $L^2_*$.\\
%Now assume $Gu=0$. Then we have:
%\[P_{\tau}u-u=\int_0^{\tau}P_sGu ds=0 \quad \mbox{for all}\; \tau>0  \]
 %Hence, the kernel of $G$ is one-dimensional and contains exactly the constants. 
 %Let now $u \in \mbox{Ker} G^2$ that is $Gu\in \mbox{Ker}G$. Thus, we must have $Gu\equiv c$ for some constant $c$.  
 %But since $\int Gu dv =0$ by \ref{diss} we can deduce $c=0$, thus $u$ already was in Ker $G$ and we have Ker $G^2=$Ker$G$.

\section{Asymptotic Behaviour of the Semigroup}\label{asy}

Having obtained a unique invariant measure for the semigroup of the reduced equation, we may use ergodic theory to ivestigate the asymptotics of the two-parameter semigroup. 
\begin{prop}
Assume that $f:H\to \R$ is such that\\
 $\int_0^T \int_H f^2(x) \nu_t(dx) dt <\infty \quad$ Then we have:
$$ \lim_{\tau \to\infty}\frac{1}{\tau}\int_0^{\tau} P_{t,t+s}f(x) ds =  \frac{1}{T} \int_0^T \int_H f(x) \nu_t(dx) dt $$
\end{prop}
\begin{proof}A simple application of the ergodic theorem.
\end{proof}
%Since every unique invariant measure is automatically ergodic, the result follows directly from the equality (see \cite{erg} 3.2.4):
 %$$ \lim_{\tau \to\infty}\frac{1}{\tau}\int_0^{\tau} P_{s}g(t,x) ds =   \int_0^T \int_H g(t,x) \nu(dx) $$
 %which is valid for all $g\in L^2_*(\nu).$ We have only to take $g(t,x)=f(x)$ independent of $t$ and recall, that then $P_sg(t,x)=P_{t+s,t}f(x)$.
    
%Under the condition of weak convergence, we are able to characterize the asymptotic behaviour of $P_{s,t}$ for $s \to -\infty$ with the help of our evolution system of measures.
\begin{prop}
Assume that there is $x\in H$ such that for $s\to -\infty$\\
$P\circ (X(t,s,x))^{-1} \to \nu_t$ weakly. Then we have for $f\in C_b(H)$ :
$$\lim_{s\to -\infty} P_{s,t}f(x) = \int_H f(x)\nu_t(dx) $$ 
\end{prop}

\begin{proof}
By definition of weak convergence.
\end{proof}

\subsection{The Square Field Operator and an Estimate}\label{sqfo} 

In the following we will introduce the square field operator. Its importance lies in the crucial role, that it will play in the proof of the following functional inequalities.
 
\begin{defi}
$\Gamma(u,u):=Gu^2-2uGu$ will be called the square field operator.
\end{defi}

\begin{lem}[square field operator] \label{square}
On $K$ we have:
\begin{align*} \Gamma(u,u) = \left\langle \Sigma^*B^*(t)\nabla_x u,\Sigma^*B^*(t)\nabla_x u \right\rangle
 + \int_H \left[u(x+B(t)y,t)-u(x,t)\right]^2 M(dy)\end{align*}
\end{lem}
\begin{proof}
Let $u$ be given by $u(t,x)=\Phi(t)e^{i\langle x,h(t)\rangle}$.\\
First note, that $u^2(t,x)= \Phi^2(t)e^{i\langle x,2h(t)\rangle}$, so that by \eqref{generator}:
\begin{align}
Gu^2(t,x)&=[2\Phi'(t)\Phi(t) + i\Phi^2(t)\langle x,2h'(t)\rangle] \,e^{i\langle x,2h(t)\rangle}\label{aaa} \\
                                            &+ i\langle A(t)x + f(t),2h(t)\rangle\, \Phi^2(t)\,e^{i\langle x,2h(t)\rangle}\label{bbb}\\
                                            &+ \lambda[B^*(t)2h(t)]\, \Phi^2(t)\,e^{i\langle x,2h(t)\rangle}\label{ccc}
\end{align}
\begin{align}
2u(t,x)Gu(t,x)&=2[\Phi'(t) + i\Phi(t)\langle x,h'(t)\rangle] \,\Phi(t)\left(e^{i\langle x,h(t)\rangle}\right)^2 \label{aaaa}\\
                                            &+2 i\langle A(t)x + f(t),h(t)\rangle\, \Phi^2(t)\,\left(e^{i\langle x,h(t)\rangle}\right)^2\label{bbbb}\\
                                            &+ 2\lambda[B^*(t)h(t)]\, \Phi^2(t)\,\left(e^{i\langle x,h(t)\rangle}\right)^2\label{cccc}
\end{align}
We see immediately, that \eqref{aaa}=\eqref{aaaa} and \eqref{bbb}=\eqref{bbbb}. Then a tedious, but simple calculation shows, that \eqref{ccc} - \eqref{cccc} has indeed the proclaimed form.

Note that - since the square field operator is not linear - we have to check the equality also for sums of such functions. Another lengthy calculation yields the result.  
\end{proof}
%VIELLEICHT NOCH DISIPATIVITÄT?
%\newpage
\begin{ass}\label{hypo}\mbox{}
 \begin{itemize}
  \item[(i)] For every $t,\tau>0: U(t+\tau,t)RH\subset \sqrt{R}H$ and there is a strictly positive $C_1 \in C[0,\infty)$ such that:
  $$\|U(t,s)Rx\|_{H_0} \leq \sqrt{C_1(t-s)} \|Rx\|_{H_0} \quad x\in H,\;t>s$$
   \item[(ii)] There is a strictly positive $C_2\in C[0,\infty)$ such that:
   $$M\circ U(t+\tau,t)^{-1}\leq C_2(\tau) M \quad \tau>0$$
   that is   $C_2(\tau) M -M\circ U(t+\tau,t)^{-1}$ is a positive measure. 
\end{itemize}
\end{ass}

\begin{lem}[estimate of the square field operator] If $B=Id$, we have for $u\in K$:
\begin{equation}\label{local}
\sqrt{\langle D_xP_{\tau}u(t,x),RD_xP_{\tau}u(t,x)\rangle} \leq \sqrt{C_1(\tau)} \;P_{\tau} \left(\|\sqrt{R}D_x u\|\right)(t,x)
\end{equation}
\vspace{-5mm}
\begin{multline}\label{nonlocal}
\int_H \left[P_{\tau} u(x+y,t)-P_{\tau}u(x,t)\right]^2 M(dy) \\
\leq C_2(\tau) P_{\tau}\left(\int_H  [u(\cdot+y)-u(\cdot)]^2 M(dy)\right)(x,t)
\end{multline}
So that combining the two estimates, we have:
\[\Gamma(P_{\tau}u,P_{\tau}u)\leq \max(C_1,C_2)(\tau) P_{\tau} \Gamma(u,u)\]
\end{lem}
\begin{proof}
Let be $z\in H$ and $u(t,x)=\Phi(t)e^{i\langle x,h(t)\rangle}$, then:
\begin{align*} 
&\langle D_xP_{\tau}u(t,x),Rz\rangle = \langle i U^*(t+\tau,t)h(t+\tau)P_{\tau}u(t,x),Rz\rangle \\
%=  \left\langle i U^*(t+\tau,t)h(t+\tau)\int_H u(t+\tau,y)P\circ X(t+\tau,t,x)^{-1}(dy),Rz\right\rangle \\
 &   =\int_H \left\langle i U^*(t+\tau,t)h(t+\tau)u(t+\tau,x),Rz\right\rangle P\circ X(t+\tau,t,x)^{-1}(dy)\\ \tag{+} \label{plus}
%=P_{\tau} \langle  i U^*(t,t-\tau)h(t)u(t,x),Rz\rangle  
&=P_{\tau} \langle D_x u(t,x),\sqrt{R}\sqrt{R}^{\,-1}U(t,t-\tau)Rz\rangle \\
%=P_{\tau} \langle \sqrt(R)D_x u(t,x),\sqrt{R}^{\,-1}U(t,t-\tau)Rz\rangle \\ 
&\leq P_{\tau} \|\sqrt(R)D_x u\| \|\sqrt{R}^{\,-1}U(t,t-\tau)Rz\| 
= P_{\tau} \|\sqrt(R)D_x u\| \|U(t,t-\tau)Rz\|_{H_0} \\
&\leq P_{\tau} \|\sqrt(R)D_x u(t,x)\| \sqrt{C_1(\tau)}\|Rz\|_{H_0} 
= P_{\tau} \|\sqrt(R)D_x u(t,x)\| \sqrt{C_1(\tau)}\|\sqrt{R}z\| \\
\end{align*}
Now, for every pair $(t,x)$ choosing $z=D_xP_{\tau}u(t,x)$ we obtain:
\begin{multline*}  \langle D_xP_{\tau}u(t,x),RD_xP_{\tau}u(t,x)\rangle \\
\leq \sqrt{C_1(\tau)}\;\|\sqrt{R}D_xP_{\tau}u(t,x)\|\ \;P_{\tau} \left(\|\sqrt{R}D_x u\|\right)(t,x)\end{multline*}
or
\begin{equation}\label{gradientestimate}  \sqrt{\langle D_xP_{\tau}u(t,x),RD_xP_{\tau}u(t,x)\rangle} \leq \sqrt{C_1(\tau)} \;P_{\tau} \left(\|\sqrt{R}D_x u\|\right)(t,x)\end{equation}
Note that we have used the special form of $u$ only up to equation \eqref{plus}, but by linearity of $P$ and $D_x$ it is clear that this also holds for sums. 
So we obtain \eqref{gradientestimate} on all of $K$. %\\[5mm]

Setting $\tilde{P}:=P\circ X(t+\tau,t,0)^{-1}$ and $\tilde{M}:=M\circ U(t+\tau,t)^{-1}$ we have for general $u\in K$: (setting $\tilde{\tau}:=t+\tau$ for brevity)
\begin{align*}
 &\int_H \left|P_{\tau} u(x+y,t)-P_{\tau}u(x,t)\right|^2 M(dy) \\
%=&\int_H\left|\int_H  u(U(\tilde{\tau},t)(x+y)+z,\tilde{\tau})-u(U(\tilde{\tau},t)x+z,\tilde{\tau}) \tilde{P}(dz)\right|^2    M(dy) \\
\leq&\int_H\left(\int_H  |u(U(\tilde{\tau},t)(x+y)+z,\tilde{\tau})-u(U(\tilde{\tau},t)x+z,\tilde{\tau})|^2 \tilde{P}(dz)\right)    M(dy) \\
%=&\int_H\left(\int_H  |u(U(\tilde{\tau},t)(x+y)+z,\tilde{\tau})-u(U(\tilde{\tau},t)x+z,\tilde{\tau})|^2 M(dy)\right) \tilde{P}(dz)    \\
=&\int_H\left(\int_H  |u(U(\tilde{\tau},t)x+y+z,\tilde{\tau})-u(U(\tilde{\tau},t)x+z,\tilde{\tau})|^2 \tilde{M}(dy)\right) \tilde{P}(dz)    \\
\leq & C_2(\tau)\int_H\left(\int_H  |u(U(\tilde{\tau},t)x+y+z,\tilde{\tau})-u(U(\tilde{\tau},t)x+z,\tilde{\tau})|^2 M(dy)\right) \tilde{P}(dz)\\
=& C_2(\tau) P_{\tau}\left(\int_H  |u(\cdot+y)-u(\cdot)|^2 M(dy)\right)(x,t)
\end{align*}

\end{proof}

\begin{cor}
\begin{equation}\label{localohner}
\sqrt{\langle D_xP_{\tau}u(t,x),D_xP_{\tau}u(t,x)\rangle} \leq  \;\|U(t+\tau,t)\|P_{\tau} \left(\|D_x u\|\right)(t,x)
\end{equation}
\end{cor}
\begin{proof}
Reconsidering the proof above and setting $R=Id$ yields the result.
\end{proof}
%\newpage

\subsection{Functional Inequalities}

Following \cite{wang}, we will now prove a Poincaré and a Harnack inequality.

\begin{defi}
\[ \overline{u}_t:=\int_H u(t,x) \nu_t(dx)\]
\end{defi}
\begin{prop}\label{asympt}
Assume that the $\nu_t$ have uniformly bounded first moments, that is: 
$ \sup_t \{\int_H \|x\|\nu_t(dx)\} < \infty \quad$ Then we have for all $u\in K$: \\
\[ \lim_{\tau\to\infty}\left( \sup_t |P_{\tau}u(t,x)-\overline{u}_{t+\tau} |  \right)=0 \quad \mbox{for every fixed}\;x\]
\end{prop}
\begin{proof}
We have, since $\overline{u}_{t+\tau}:=\int_H u(t+\tau,y) \nu_{t+\tau}(dy)=P_{t,t+\tau} u(t+\tau,\cdot)(y)\nu_t(dy)$ by the property of the evolution system :
\begin{align*}&|P_{\tau}u(t,x)-\overline{u}_{t+\tau} |=\left|\int_H P_{t,t+\tau} u(t+\tau,\cdot)(x)-P_{t,t+\tau} u(t+\tau,\cdot)(y)\nu_t(dy) \right|\\
& \leq \|D_xP_{t,t+\tau} u(t+\tau,\cdot)\|_{\infty} \int_H |x-y|\nu_t(dy) 
%\\&\leq \|U(t+\tau,t)\| \|P_{\tau}D_x u(t,\cdot)\|_{\infty} \int_H |x-y|\nu_t(dy)
\\ & \leq Me^{-\omega \tau} \|D_x u(t+\tau,\cdot)\|_{\infty} \int_H |x-y|\nu_t(dy) \stackrel{\tau\to\infty}{\longrightarrow} 0 
\end{align*}
since the integral is bounded by assumption and $\|D_xu(t,x)\|= \|h(t)\| $ but $h: \R\to H$ is continuous and periodic, hence bounded.  
\end{proof}

\begin{prop}[Poincaré Inequality]\label{poin}
Given assumption \ref{hypo} \\ and $B=Id$, we have for $C(\tau):=\max\left(\int_0^{\tau} C_1(s)ds,\int_0^{\tau} C_2(s)ds\right):$
\begin{equation} \label{poineq} P_{\tau}u^2-(P_{\tau}u)^2 \leq C(\tau)P_{\tau}\Gamma(u,u)  \qquad \mbox{for all}\quad\tau>0 , u\in K  \end{equation}
\end{prop}
\begin{proof}
Set $f(s):=P_{\tau-s}(P_{s}u)^2$ Then we have by the product rule:
\begin{align*} \frac{d}{ds} f(s)&= -P_{\tau-s}G(P_{s}u)^2 +P_{\tau-s}2P_{s}uGP_su \\ &=-P_{\tau-s}[G(P_{s}u)^2-2P_{s}uGP_su]=-P_{\tau-s}\Gamma(P_su,P_su)\end{align*}
Hence,
\begin{align*}
-\frac{d}{ds} f(s) & = P_{\tau-s}\Gamma(P_su,P_su)&\\
                   & = P_{\tau-s} \left\langle \Sigma^*\nabla_x P_su,\Sigma^*\nabla_x P_su \right\rangle &\\
                   &\quad + P_{\tau-s}\int_H \left[P_su(x+y,t)-P_su(x,t)\right]^2 M(dy)  &\\
                   & \leq  C_1(s)  P_{\tau-s}P_s \left\langle \Sigma^*\nabla_x u,\Sigma^*\nabla_x u \right\rangle & \mbox{by \eqref{local}}\\
                   & \quad + C_2(s) P_{\tau-s}P_s\int_H \left[u(x+y,t)-u(x,t)\right]^2 M(dy) &  \mbox{by \eqref{nonlocal}}\\
\end{align*}
Integrating with respect to $s$ and noting that $f(0)=P_{\tau}f^2$ and $f(t)=(P_{\tau}f)^2$ we obtain:
\begin{align*} P_{\tau}f^2-(P_{\tau}f)^2 &\leq \left(\int_0^{\tau} C_1(s)  ds\right) P_{\tau} \left\langle \Sigma^*\nabla_x u,\Sigma^*\nabla_x u \right\rangle \\
& \quad +  \left(\int_0^{\tau} C_2(s) ds\right)  P_{\tau}\int_H \left[u(x+y,t)-u(x,t)\right]^2 M(dy) \end{align*}
and the result is proved.
\end{proof}

\begin{cor}Let be $C$ as in \ref{poin}.\\
Given that $C(\infty)<\infty$ we also have for all $\,u \in K$:
\[\int_{[0,T]\times H} [u(t,x)-\overline{u}_t]^2\nu(dt,dx) \leq C(\infty) \int_{[0,T]\times H} \Gamma(u,u) \nu(dt,dx)\]  
\end{cor}
\begin{proof}
Integrating \eqref{poineq} with respect to $\nu$ yields, because of invariance:
\[\int_{[0,T]\times H}u^2-(P_{\tau}u)^2\nu(dt,dx) \leq  C(\tau) \int_{[0,T]\times H} \Gamma(u,u) \nu(dt,dx) \]
Letting $\tau \to\infty$ and using \ref{asympt} together with dominated convergence we have:
\[\int_{[0,T]\times H}u^2-(\overline{u}_{t})^2\nu(dt,dx) \leq  C(\infty) \int_{[0,T]\times H} \Gamma(u,u) \nu(dt,dx) \]
Since $\overline{u}$ does not depend on $x$ anymore, we have:
\begin{align*}&\int_{[0,T]\times H} [u(t,x)-\overline{u}_t]^2\nu(dt,dx)=\int_{[0,T]\times H} u^2(t,x)-2u(t,x)\overline{u}_t+\overline{u}^2_t\nu(dt,dx)\\
&=\int_{[0,T]\times H} u^2(t,x)+\overline{u}^2_t\nu(dt,dx) -2\int_{[0,T]}\overline{u}_t\underbrace{\int_H u(t,x) \nu_t(dx)}_{\overline{u}_t}dt
%\\&=\int_{[0,T]\times H} u^2(t,x)-\overline{u}^2_t\nu(dt,dx)
\end{align*}
%and the result follows.
\end{proof}
%\newpage
 For the following Harnack inequality we need a definition:
 \begin{defi}
 $$\rho(x,y):=\inf \{\; \|z\|\;:\;  \sqrt{R}z=x-y \}$$
 with the usual convention that $\;\inf \emptyset = \infty\;$, so $\rho$ may take the value infinity if $(x-y)\notin$ Im$\sqrt{R}$.\\
 
 $C_b^*:=\{f:\R \times H \to \R \;|\, f \mbox{is continuous and bounded}$\\
 $\mbox{}\qquad\qquad\qquad\qquad\qquad\qquad\mbox{ and T-periodic in the first component}\}$
 \end{defi}
            
\begin{prop}[Harnack Inequality]
\begin{equation}\label{harnack}|P_{\tau} u(t,y)|^2 \leq P_{\tau}u^2(t,x) \exp\left[ \frac{\rho^2(x,y)}{\int_0^{\tau} \frac{1}{h(s)}ds}\right] \quad \mbox{for all} \;u\in C_b^* \end{equation}
\end{prop}
\begin{proof}
First, let be $u\in K$ such that $u$ is strictly positive.
Since \\ $P_{\tau-s}(P_{s}u)^2(t,x)$ will then also be strictly positive we can define:
$$ \Phi(s):= \log[P_{\tau-s}(P_{s}u)^2(t,x_s)]$$
where $x_s$ is given by
$ x_s:= x+\frac{(y-x)\int_0^s \frac{1}{h(\tau-u)} du}{\int_0^{\tau} \frac{1}{h(u)}du} \quad$
%Note that we have $x_0=x$ and $x_{\tau}=y$.\\
Differentiating $\Phi$ we obtain:
\begin{align}\label{h1}
\frac{d}{ds}\Phi(s)&=\frac{\frac{d}{ds}P_{\tau-s}(P_{s}u)^2(t,x_s)}{P_{\tau-s}(P_{s}u)^2(t,x_s)}
\end{align}
and for the numerator:
\begin{align}\label{h2}
\frac{d}{ds}[P_{\tau-s}(P_{s}u)^2(t,x_s)]&=\frac{d}{ds}[P_{\tau-s}(P_{s}u)^2](t,x_s) + \left\langle D_x[P_{\tau-s}(P_{s}u)^2](t,x_s), \frac{dx_s}{ds}\right\rangle \notag \\
                                         &= -GP_{\tau-s}(P_su)^2(t,x_s) + P_{\tau-s}[2P_suGP_su](t,x_s)  \notag\\
                                         &\quad + \frac{1}{h(\tau-s)\int_0^{\tau} \frac{1}{h(u)} du}\left\langle  D_x[P_{\tau-s}(P_{s}u)^2](t,x_s),(y-x)\right\rangle  \notag \\  
&=  -P_{\tau-s}\Gamma(P_{s}u,P_{s}u)\notag\\
&\quad + \frac{1}{h(\tau-s)\int_0^{\tau} \frac{1}{h(u)} du}\left\langle  D_x[P_{\tau-s}(P_{s}u)^2](t,x_s),(y-x)\right\rangle \notag \\
\end{align}
We will now estimate $\left\langle  D_x[P_{\tau-s}(P_{s}u)^2](t,x_s),(y-x)\right\rangle$:

\begin{align}\label{h3}
&\left\langle  D_x[P_{\tau-s}(P_{s}u)^2](t,x_s),(y-x)\right\rangle  &\notag\\
&=\inf_{\{z: \sqrt{R}z=x-y\}}\left\langle  D_x[P_{\tau-s}(P_{s}u)^2](t,x_s),\sqrt{R}z\right\rangle  \qquad \qquad \qquad (x-y)\in \mbox{Im}\sqrt{R} \notag \\
&\leq \sqrt{\left\langle  RD_x[P_{\tau-s}(P_{s}u)^2](t,x_s),D_x[P_{\tau-s}(P_{s}u)^2](t,x_s)\right\rangle}\rho(x,y)  \qquad \; \mbox{Cau.-Schw.}\notag\\
&\leq \rho(x,y) \sqrt{h(\tau-s)}P_{\tau-s}\left(\sqrt{\left\langle  RD_x(P_{s}u)^2,D_x(P_{s}u)^2\right\rangle}\right)(t,x_s) \;\;\quad \quad \mbox{by \eqref{local}} \notag\\
&\leq 2\rho(x,y)\sqrt{h(\tau-s)} P_{\tau-s}\left(P_{s}u\sqrt{\left\langle  RD_x(P_{s}u),D_x(P_{s}u)\right\rangle}\right)(t,x_s)  \quad \mbox{chain rule}  \notag\\ 
\end{align}
Combining \eqref{h1},\eqref{h2} and \eqref{h3} we obtain:

\begin{align*}
&\frac{d}{ds}\Phi(s)
\leq  \qquad \frac{-P_{\tau-s}\Gamma(P_{s}u,P_{s}u)}{P_{\tau-s}(P_{s}u)^2(t,x_s)}\\[0mm]
                   &+\frac{\frac{1}{h(\tau-s)\int_0^{\tau} \frac{1}{h(u)} du} 2\rho(x,y)\sqrt{h(\tau-s)} P_{\tau-s}\left(P_{s}u\sqrt{\left\langle  RD_x(P_{s}u),D_x(P_{s}u)\right\rangle}\right)(t,x_s) }
{P_{\tau-s}(P_{s}u)^2(t,x_s)} \\[0mm]
%& \mbox{\textit{writing out}} \;\Gamma\; \mbox{\textit{but omitting the non-local part, we get:}}\\[5mm]
                   &\leq \frac{-P_{\tau-s}\left(\left\langle  RD_x(P_{s}u),D_x(P_{s}u)\right\rangle\right)(t,x_s)}{P_{\tau-s}(P_{s}u)^2(t,x_s)}\\[0mm]
                   &+\frac{\frac{1}{\sqrt{h(\tau-s)}\int_0^{\tau} \frac{1}{h(u)} du} 2\rho(x,y) P_{\tau-s}\left(P_{s}u\sqrt{\left\langle  RD_x(P_{s}u),D_x(P_{s}u)\right\rangle}\right)(t,x_s) }
{P_{\tau-s}(P_{s}u)^2(t,x_s)}\\[0mm]
                   &= \frac{1}{P_{\tau-s}(P_{s}u)^2(t,x_s)}\\
                   &\times P_{\tau-s}\!\left((P_{s}u)^2 \! \left[
                     2H\frac{\sqrt{\left\langle  RD_x(P_{s}u),D_x(P_{s}u)\right\rangle}}{P_{s}u}-\frac{\left\langle  RD_x(P_{s}u),D_x(P_{s}u)\right\rangle}{(P_{s}u)^2} \right]\right)\!(t,x_s)\\
\end{align*}
where we have set $H:=\frac{\rho(x,y)}{\sqrt{h(\tau-s)}\int_0^{\tau} \frac{1}{h(u)} du} $ for brevity.\\[3mm]
Furthermore, setting $G:=\frac{\sqrt{\left\langle  RD_x(P_{s}u),D_x(P_{s}u)\right\rangle}}{P_{s}u}$ :
\begin{align*} \frac{d}{ds}\Phi(s)&  \leq \frac{1}{P_{\tau-s}(P_{s}u)^2(t,x_s)}
                                     P_{\tau-s}\left((P_{s}u)^2\left[ -G^2+ 2HG \right]\right)(t,x_s)\\
                                  & = \frac{1}{P_{\tau-s}(P_{s}u)^2(t,x_s)}
                                     P_{\tau-s}\left((P_{s}u)^2\left[ -G^2+ 2HG -H^2+H^2\right]\right)(t,x_s)\\
                                  & \leq \frac{1}{P_{\tau-s}(P_{s}u)^2(t,x_s)}
                                     P_{\tau-s}\left((P_{s}u)^2\left[ H^2 \right]\right)(t,x_s)\\
                                     &= H^2
\end{align*}
since $H$ depends neither on $x_s$ nor on $t$. Integration over $s$ yields:
\begin{align*}
&\log[ (P_{\tau}u)^2(t,y)]-\log[(P_{\tau}u^2)(t,x)]=\Phi(\tau)-\Phi(0) \\[2mm]
&\leq \int_0^{\tau} H^2(s)ds =\int_0^{\tau}\frac{\rho^2(x,y)}{h(\tau-s)(\int_0^{\tau} \frac{1}{h(u)} du)^2} ds =\frac{\rho^2(x,y)}{\int_0^{\tau} \frac{1}{h(u)} du} 
\end{align*}
Hence, applying the exponential yields:
$(P_{\tau}u)^2(t,y)  \leq  P_{\tau}u^2(t,y) \; \frac{\rho^2(x,y)}{\int_0^{\tau} \frac{1}{h(u)} du} $\\
and the proof is complete for positive functions.
To obtain the result for general $u$, note first, that it is sufficient to have it for $|u|$, since we have:\\
$|P_{\tau} u(t,y)|^2 \leq [P_{\tau} |u|(t,y)]^2 \leq P_{\tau}u^2(t,x) \exp\left[ \rho^2(x,y)(\int_0^{\tau} \frac{1}{h(s)}ds)^{-1}\right]$\\
Of course, we cannot take modulus without leaving $K$, but as $K$ is an algebra we may take the square of our functions. Thus, let be $u\in C_b^*$ and $\eps>0$. Then $f:=\sqrt{|u|}\in C_b^*$. Now, by \ref{approx} we can approximate $f$ pointwisely by functions $u_n$ from $K$. Then, $u_n^2 +\eps $ is strictly positive, it will approach $|u|+\eps$ and since the approximating functions are uniformly bounded, we can take limits in \eqref{harnack} and obtain the result via dominated convergence and then letting $\epsilon\to 0$.   
\end{proof}

\begin{lem}\label{approx}
For every $f\in C_b^*$ we can find a sequence $u_n \in K$ such that:
\begin{itemize}
\item $u_n \to f $ pointwisely
\item $ \sup_{x,t,n} |u_n(t,x)| \leq 1+\sup_{t,x} |f(t,x)|$ 
\end{itemize}
\end{lem}
\begin{proof}
Let be $(e_k)_{k\in\N}$ a complete orthonormal system in $D(A^*)$.
Let be \\
$f_n(t,h):=g_n(t,P_nh)$
 where $P_n: H\to \R^n \quad h\mapsto (\langle h,e_1\rangle, \dots ,\langle h,e_n\rangle)$\\
 and $g_n:\R \times \R^n\to \R\quad (t,x_1,\dots,x_n)\mapsto f(t,x_1e_1+\dots+x_ne_n)$
 
 Note that each $g_n$ is continuous and bounded and that $f_n=f$ on \\ $\R \times \mbox{span}\{e_1,\dots,e_n\}$.
 Moreover, we have $P_n\to Id $ strongly and hence $f_n\to f$ pointwisely.
 
 Now let be $g:\R \times \R^d\to \R$ continuous and bounded. 
 By the theorem of Stone-Weierstrass it is clear, that we can approximate $g$ uniformly on any bounded set by sums of products of complex exponentials with T-periodic differentiable functions.

 %For each $n\in \N$ let be be $\phi_n$ a function that coincides with $g$ 
 %on $\{\|x\|_{\infty}\leq n\}$ vanishes on $\{\|x\|_{\infty}\geq n+1\}$ and is still continuous.
 %Since $\phi_n$ can be extended to a continuous $(n+1)$-periodic function,  we can approximate it with respect to the supremum norm on $\{\|x\|_{\infty}\leq n+1\}$  with trigonometric polynomials, by Fejér's Theorem.  
 Let us call $g^{(n)}$ such an approximation of $g$ on $\{\|x\|_{\infty}\leq n\}$ up to $\frac{1}{n}$, that is we have:
 $\sup_{\|(t,h)\|_{\infty}\leq n}|g^{(n)}(t,h)-g(t,h)|<\frac{1}{n}$  so we have $g^{(n)} \to g$ pointwisely for $n\to \infty$ and by periodicity it is clear that we have $ \sup_{t,x,n} |g^{(n)}(t,x)| \leq 1+\sup_{t,x} |g(t,x)|$.
 
 Now, applying the above approximation to our functions $g_n$ from above, we denote: 
 \[u_n(t,h):=f_n^{(n)}(t,h):=g_n^{(n)}(t,P_nh)\]
 and an easy calculation shows that $f_n^{(n)}$ is indeed a function from $K$. 
 %, since we have:
 %\begin{align*}e^{i\langle (k_1,\dots,k_n),P_nh\rangle}&=e^{i\langle (k_1,\dots,k_n),(\langle h,e_1\rangle,\dots,\langle h,e_n\rangle)\rangle}\\
 %&= e^{i(\langle h,k_1e_1\rangle+\dots+\langle h,k_ne_n\rangle)}=e^{i\langle h,k_1e_1+\dots+k_ne_n\rangle}\end{align*}
%and  $(k_1e_1+\dots+k_ne_n)\in D(A^*)$.\\
 To see that $f_n^{(n)}(t,h) \to f(t,h)$ for each fixed $(t,h)$ we calculate:
 \[ |f(t,h)-f_n^{(n)}(t,h)|\leq   |f(t,h)-f_n(t,h)|+ |f_n(t,h)-f_n^{(n)}(t,h)|\]
We have already stated that the first tern will tend to $0$. For the second term, note, that for fixed $(t,h)$ there is $N$ independent of $n$ such that 
$ \|(t,P_n h)\|_{\infty}<N$. Thus, we obtain:
 \begin{gather*}|f_n(t,h)-f_n^{(n)}(t,h)|=|g_n(t,P_nh)-g_n^{(n)}(t,P_nh)|\\
 \leq \sup_{\|(t,x)\|_{\infty}\leq N}|g_n(t,x)-g_n^{(n)}(t,x)|<\frac{1}{n}\end{gather*}
 whenever $n$ is bigger than $N$.

\end{proof}

%\newpage
\section{Appendix}
\vspace{2mm}

\hspace{-2mm}\textbf{\large{ Proof of Lemma \ref{ft}:}}\\[2mm]
 Knowing how the Fourier transform acts on translations, it will be enough to show, that:
$$\E\left[\exp\left(i\left\langle h,\int_s^t U(t,r)B(r) dL_r\right\rangle \right)\right]= \exp \left\{ \int_s^t \lambda(B^*(r)U^*(t,r)h)dr \right\}$$
The strong continuity of $U$ and $B$ allows us to approximate the Lévy stochastic integral by a sequence of sums. More precisely, we have:
\[\int_s^t U(t,r)B(r) dL_r = P-\limes \sum_{s_i\in \mathcal P_n} U(t,s_i)B(s_i)(L_{s_i}-L_{s_{(i-1)\vee 0}})\]
where the limit is taken in probability and $\mathcal P_n$ is a sequence of partitions $s=s_0<...<s_N=t$ of $[s,t]$ such that the mesh width tends to zero.
The verification is a straightforward application of the respective isometries: 

For the drift term which is a Bochner integral we have to show that:
\[\limes \sum_{s_i \in \mathcal P_n} \int_{s_{i-1}}^{s_i} \| U(t,s_i)B(s_i)b-U(t,r)B(r)b  \|\, dr =0\]
but since $r\mapsto U_t(r) B(r)b$ is even uniformly continuous on $[s,t]$ we may find $\delta >0$ such that 
 $ \|U_t(r) B(r)b -U_t(r') B(r')b \| < \frac{\epsilon}{t-s} $ whenever $|r-r'|<\delta$, so that if we choose $n$ such that the mesh width of $\mathcal P_n$ is smaller than $\delta$ we have 
 $$\sum_{s_i \in \mathcal P_n} \int_{s_{i-1}}^{s_i} \| U(t,s_i)B(s_i)b-U(t,r)B(r)b  \|\, dr <\sum_{s_i \in \mathcal P_n} \int_{s_{i-1}}^{s_i} \frac{\epsilon}{t-s} dr < \epsilon $$

 For the small jumps we make use of the isometry from \ref{isometry}, so we have to show that our piecewise approximation converges in the $\mathcal H^2$ norm, that is we need:
 \[ \limes \sum_{k=1}^{\infty}\int_{\|x\|<1}\left(\sum_{s_i\in\mathcal P_n} \int_{s_{i-1}}^{s_i} \| [U_t(r)B(r)-U_t(s_i)B(s_i)] T_x^{\frac{1}{2}}e_k   \|^2 dr \right) M(dx)  = 0 \]
 For each $k$ and $x$ fixed the expression in round brackets converges to zero, for the same reasons as used for the drift term. So we only have to show that we may take the limit into the sum and the integral, but this follows by dominated convergence on considering the uniform integrable bound :
 \[\| [U_t(r)B(r)-U_t(s_i)B(s_i)] T_x^{\frac{1}{2}}e_k \|^2   \leq   2 \sup_{s\leq r\leq t} \|U_t(r)\|_{\mathcal L(H)} \sup_{s\leq r\leq t} \|B(r)\|_{\mathcal L(H)} \|T_x^{\frac{1}{2}}e_k\|^2   \]
 Thus we have convergence in $L^2$ of the approximating sums towards the integral.\\
 The same argument works for the Brownian part, where there is even no dependence on $x$.\\
 The big jumps, finally are quite simple to treat. Since the expression makes sense pointwise, we consider the approximation for $\omega$ fixed and we obtain: 
 \[ \limes      \sum_{s_i\in\mathcal P_n} \sum_{s_{i-1}\leq r \leq s_i} [U_t(s_i)B(s_i)-U_t(r)B(r)] \Delta L_r(\omega) \chi_{\|\Delta L_r(\omega)\|>1}=0  \]  
 again because of strong continuity.
 
 So in any of the four cases we have  at least convergence in probability and the claim is proved.

By convergence in distribution and independence of increments:
\begin{gather*}  \E\left[\exp\left(i\left\langle h,\int_s^t U(t,r)B(r) dL_r\right\rangle \right)\right]\\
%\[= \lim_{n \to \infty}\E\left[  \exp\left(i\left\langle h,\sum_{k\in\mathcal P_n} U(t,s_k)B(s_k)(L_{s_k}-L_{s_{k-1}\vee 0})\right\rangle \right)  \right]\]
 =\lim_{n \to \infty}\prod_{k\in\mathcal P_n}\E\left[  \exp\left(i\left\langle h, U(t,s_k)B(s_k)(L_{s_k}-L_{s_{k-1}\vee 0})\right\rangle \right)  \right]\\
 %=\lim_{n \to \infty} \prod_{k\in\mathcal P_n} \exp \left\{ \lambda(B^*(s_k)U^*(t,s_k)h) (s_k-(s_{k-1}\vee 0))\right\} \\
 = \exp \left\{ \int_s^t \lambda(B^*(r)U^*(t,r)h)dr \right\}  \end{gather*}
 where we employed the Lévy-Khinchine formula and the functional equation of the exponential. Note that the Riemannian sums converge to the integral because of strong continuity.\\[5mm]
 \hspace{-2mm}\textbf{\large{ Proof of Lemma \ref{fre}:}}\\[2mm]
%\begin{Proof}[]:
Let be $\lambda$  the corresponding Lévy symbol and $M$ the Lévy measure.
By the Lévy-Khinchine formula  \ref{kin} we know that:
\[ \lambda (u)= i\langle u,b \rangle -\frac{1}{2} \langle u,Au \rangle +
\int\left(e^{i\langle u,x \rangle} -1 -i \langle u,x \rangle \chi_{\{\|x\|\leq 1\}} \right)M(dx) \]
Clearly, it is enough to show that the integral expression is differentiable. We first show Gâteaux differentiability, hence we will need the directional derivatives to be integrable to obtain the result via dominated convergence. We have:
\[\left.\frac{\partial}{\partial t}\right|_{t=0} \left(e^{i\langle u+tv,x \rangle} -1 -i \langle u+tv,x \rangle \chi_{\{\|x\|\leq 1\}} \right) 
= i\langle v,x\rangle e^{i\langle u,x \rangle}  -i\langle v,x\rangle \chi_{\{\|x\|\leq 1\}} \]
To see the integrability we split the integral in two parts:
\begin{align*}
&\int_{\|x\|\leq 1} \left| i\langle v,x\rangle e^{i\langle u,x \rangle}  -i\langle v,x\rangle \chi_{\{\|x\|\leq 1\}}\right| M(dx)\\
&=\int_{\|x\|\leq 1}  \left|i\langle v,x\rangle \sum_{k=0}^{\infty} \frac{(i\langle u,x \rangle)^k}{k!} -i\langle v,x\rangle\right| M(dx)\\
%&\leq \int_{\|x\|\leq 1}  \left( \|v\|\,\|x\| \sum_{k=1}^{\infty} \frac{|(i\langle u,x \rangle)|^k}{k!} \right) M(dx)\\
%&\leq \int_{\|x\|\leq 1}  \left( \|v\|\,\|x\| \sum_{k=1}^{\infty} \frac{ \| u\|^k\|x \|^k}{k!} \right) M(dx)\\
&\leq \int_{\|x\|\leq 1}  \left(\|v\|\, \|x\|^2\|u\| \sum_{k=1}^{\infty} \frac{ \| u\|^{k-1}\|x \|^{k-1}}{k!} \right) M(dx)\\
&\leq \sup_{\|x\|\leq 1}\exp\{\|u\|\|x\|\}  \int_{\|x\|\leq 1}  \left(\|v\|\, \|x\|^2\|u\|\right) M(dx)\\
&=\exp\{\|u\|\}\|u\|\,\|v\|\, \int_{\|x\|\leq 1}  \|x\|^2 M(dx) < \infty   
\end{align*}
for every fixed $u$,\,$v$ and $s$\; since $M$ is a Lévy measure.
On the other hand, we have: 

\[ \int_{\|x\|> 1} \left| i\langle v,x\rangle e^{i\langle u,x \rangle}  -i\langle v,x\rangle \chi_{\{\|x\|\leq 1\}} \right| M(dx)
\leq \|v\|\,\int_{\|x\|> 1} \|x\| M(dx) < \infty \]
 by assumption.\\ 
Moreover, from the above, it is easy to see that the Gâteaux derivative is linear and bounded and depends continuously on $u$ with respect to the operator norm, so $\lambda$ is Fréchet differentiable and hence locally Lipschitz.\\[5mm]
%\end{Proof}
\hspace{-2mm}\textbf{\large{Proof of Proposition \ref{Kdense}}}\\[2mm]
Note that by periodicity, we can think of our functions to be defined on $[0,T]\times H$
    and in the following we will do so without changing notation.\\ 
    We will show density of $\M$ in $L^2_*(\nu;\C)$. This implies density of the respective real vector spaces.    
    We will use complex monotone classes again.
    The space $\M$ is closed under multiplication and conjugation.
    Consider $\mathcal{H}:= \bar{\M}$ as a subspace of $L^2_*(\nu;\C)$ where we allow complex-valued integrable and periodic functions.  $\mathcal{H} $ is a complex monotone vector space, by monotone convergence, applied separately to real and imaginary parts. So, $\mathcal{H}$ contains all $\sigma(\M)$-measurable functons. If we can show that $\sigma(\M)= \B(H\times\R)$, then we will have all step functions in $\mathcal H$, so density will be obvious. Note that we want to show, that functions of the form $\Phi_i \otimes e_h$ generate a product $\sigma$-algebra $\B([0,T]) \otimes \B(H)$. Since both families contain the constant function, we can break the problem down,
    as $ (1_{\R}\otimes f)^{-1}(A)=\R \times f^{-1}(A)$
     and knowing that $\Phi_i$ generates $\B(R)$ and $e_h:=e^{i\langle \,\cdot\,,h\rangle}$ generates $\B(H)$ (which follows again from \ref{borel} and the fact that $D(A^*)$ is dense), we have the result.       \\[5mm]
\newpage 
\noindent
\textbf{\large{Proof of Proposition \ref{invext}}}\\[2mm]
For invariance we have to show:
   \[
   \int_{[0,T]\times H} P_{\tau}u(t,x) d\nu =\int_{[0,T]\times H} u(t,x)d\nu \quad \forall \tau>0,\, u \in K
   \]
    Writing $ u_t(x):=u(t,x)$, remember, that $(P_{\tau} u)(t,x)= (P_{t,t+\tau}u_{t+\tau})(x)$. Taking into account \eqref{eq:esm},which is valid, since the elements from $K$ are bounded, we have:
    \[
 \int_{[0,T]\times H} P_{\tau}u(t,x) d\nu=\frac{1}{T}\int_{[0,T]}\int_{H} (P_{t,t+\tau}u_{t+\tau})(x) \nu_t(dx) dt\]
\[=   \frac{1}{T}\int_{[0,T]}\int_{H} u_{t+\tau}(x) \nu_{t+\tau}(dx) dt=\frac{1}{T}\int_{[\tau,T+\tau]}\int_{H} u_{t}(x) \nu_{t}(dx)\]
\[=    \frac{1}{T}\int_{[0,T]}\int_{H} u_{t}(x) \nu_{t}(dx)=\int_{[0,T]\times H} u(t,x)d\nu\]
because of translation invariance of $dt$ and $T$-periodicity of $u$ and $\nu_t$.

For the contraction property we have to show:  $ \|P_{\tau}u\|_{L^2_*}\leq \|u\|_{L^2_*}$\\
Using the Jensen inequality for the expectation and afterwards the invariance property for $u^2$ (recall that $K$ is closed under multiplication): 
\[ \|P_{\tau}u\|_{L^2_*}=\int_{[0,T]\times H} \E[u(t+\tau,X(t+\tau,t,x))]^2\nu(dx,dt)\]
\[ \leq \int_{[0,T]\times H} \E[u^2(t+\tau,X(t+\tau,t,x))]\nu(dx,dt) =\int_{[0,T]\times H} (P_{\tau} u^2)(t,x)\nu(dx,dt)\]
\[=\int_{[0,T]\times H} u^2(t,x)\nu(dx,dt)=\|u\|_{L^2_*}\]

To show uniqueness, let $\mu$ be another invariant measure for $P_{\tau}$, so that we have: 
\begin{equation}\label{inv}
 \int_{[0,T]\times H} P_{\tau}u(t,x) \mu(dx,dt) =\int_{[0,T]\times H} u(t,x)\mu(dx,dt) \quad \forall \tau>0,\, u \in L^2_*(\nu) \end{equation}
By \cite{dude} : corollary 10.2.8 , we can disintegrate $\mu$ as follows:
\begin{equation}\label{dis}  
\int u(t,x) \mu(dt,dx) = \int_{[0,T]} \left(\int_{H} u(t,x) \mu_t(dx)\right) \mu_1(dt)
\end{equation}
 for the marginal $\mu_1(dt)=\mu \circ Pr^{-1}$ where $Pr$ is the Projection on the $t$-component, and $\{\mu_t\}_{t\in \R}$ is a family of probability measures on $H$.
Choosing $u(t,x)=f(t)$  independent of $x$ in \eqref{inv}  we have by \eqref{dis}:
	\[\int_{[0,T]\times H} f(t+\tau) \mu_1(dt) =\int_{[0,T]\times H} f(t) \mu_1(dt)\]
	Since $f$ is $T$-periodic, $\mu_1$ is translation invariant (note, that we need here a similar monotone class argument as in \ref{Kdense}). So $\mu_1$ must be Lebesgue measure.
	
	To show $\mu_t = \nu_t$, we will of course use the uniqueness property from \ref{central}.
	Choosing $u(t,x)=f(t)g(x)$ and $\tau=T$ in \eqref{inv} yields:
	  \[\int_{[0,T]}f(t)\left(\int_{H}P_{t,t+T} g(x) \mu_t(dx)\right) \mu_1(dt) =\int_{[0,T]}f(t)\left(\int_{H} g(x)\mu_t(dx)\right) \mu_1(dt)\]
Clearly, if this holds for a fixed, bounded $g$ and arbitrary bounded $f$, we must have 
\[\int_{H}P_{t,t+T} g(x) \mu_t(dx) = \int_{H} g(x)\mu_t(dx)\]
Since this holds for any bounded measurable $g$ we can apply \ref{central} \\
to obtain $\nu_t=\mu_t \;\forall \, t \in \R$. \\[5mm]
%
%\hspace{-2mm}\textbf{\large{Proof of Proposition \ref{approx}}}\\[2mm]

\end{document}